\documentclass[a4paper,12pt]{article} 

\usepackage{color}

\usepackage{amsfonts, amsmath, amsthm, amssymb}
\usepackage[T1]{fontenc}
\usepackage[cp1250]{inputenc}

\usepackage{graphicx}
\usepackage{amssymb}
\usepackage{amsmath}
\usepackage{mathptmx}
\usepackage{helvet}
\usepackage{courier}
\usepackage{txfonts}
\usepackage{tikz} 
\usetikzlibrary{arrows}
\usepackage{type1cm}

\usepackage{verbatim}

\usepackage{graphicx}
\usepackage{epsfig,amscd,amssymb,amsxtra,amsmath,amsthm}
\usepackage{type1cm}
\usepackage[T1]{fontenc}
\usepackage{graphics}
\usepackage[mathscr]{eucal}
\usepackage[all]{xy}
\usepackage{amsmath,amscd}

\newtheorem{theorem}{Theorem}[section]
\newtheorem{proposition}[theorem]{Proposition}
\newtheorem{definition}[theorem]{Definition}

\newtheorem{example}[theorem]{Example}
\newtheorem{remark}[theorem]{Remark}

\newtheorem{corollary}[theorem]{Corollary}
\newtheorem{problem}[theorem]{Problem}
\newtheorem{observation}[theorem]{Observation}
\newtheorem{construction}[theorem]{Construction}

\newcommand{\diam} {\mathop{\rm diam}\nolimits}

\newcommand{\Cl}  {\mathop{\rm Cl}\nolimits}

\newcommand{\of}[1]{\left(#1\right)}
\newcommand{\ofb}[1]{\left[#1\right]}
\newcommand{\ofc}[1]{\left\{#1\right\}}
\newcommand{\ofi}[1]{\left<#1\right>}
\newcommand{\ofa}[1]{\left|#1\right|}

\begin{document}

\title{Characterizations of $\mathcal P$-like continua that do not have the fixed point property}
\author{Iztok Bani\v c, Judy Kennedy and Piotr Minc}
\date{}

\maketitle

\begin{abstract}
We give two characterizations of $\mathcal P$-like continua $X$ that do not have the fixed point property. Both characterizations are stated in terms of sequences of open covers of $X$ that follow fixed-point-free patterns. We use these to characterize planar tree-like continua that do not have the fixed point property in terms of infinite sequences of tree-chains in the plane that follow fixed-point-free patterns. We also establish a useful relationship between these tree-chains and commutative simplicial diagrams that we use later to construct a finite sequence (of any given length) of tree-chains in the plane that follows a fixed-point-free pattern.

An earlier characterization of $\mathcal P$-like continua with the fixed point property was given in 1994 by Feuerbacher based on a 1963 result by Mioduszewski. The Mioduszewski-Feuerbacher characterization is expressed in terms of almost commutative inverse diagrams. In contrast, our approach is more geometric, and it may potentially lead to new methods in the elusive search for a planar tree-like continuum without the fixed-point property. 
\end{abstract}
\-


\section{Introduction}\label{s1}

All spaces in this paper are metric.
Let $X$ be a continuum and let $\mathcal P$ be a class of polyhedra\footnote{The definition of a polyhedron may be found in \cite[p. 470--473]{poli}}. A finite family $\mathcal{U}$ of open sets in  $X$ is called a $\mathcal P$-{\em cover\/} for $X$ if $\bigcup_{U\in\mathcal U}U=X$ and the geometric realization of the nerve of $\mathcal{U}$ is homeomorphic to a member of $\mathcal P$. $X$ is said to be $\mathcal P$-like if and only if each open cover of $X$ can be refined by a finite open $\mathcal P$-cover of $X$.

If $\mathcal P$ is the collection of all trees, and a continuum $X$ is $\mathcal P$-like, we simply say that $X$ is tree-like.

If $X$ is $\ofc{P}$-like for some polyhedron $P$, we say that  $X$ is $P$-like. It is well known that every plane continuum $X$ does not separate the plane if and only if it is $I^2$-like. (By $I^n$ we understand the $n$ dimensional cube $\ofb{0,1}^n$.)

This paper is motivated by the plane fixed point problem: \\\centerline{``Does every nonseparating plane continuum have the fixed-point property?''}
According to Rogers \cite[p. 305]{6}, the paper of Ayres \cite{9} in 1930 was the first instance in which this problem appeared in print.
Ayres proved that each homeomorphism of a nonseparating locally connected plane continuum has a fixed point. He wrote in \cite[p. 336]{9} ``this result is a partial solution of the well-known problem as to whether a general bounded continuum not separating its plane has this property.''

\begin{figure}[h]
\centering
\scalebox{0.5}{\includegraphics{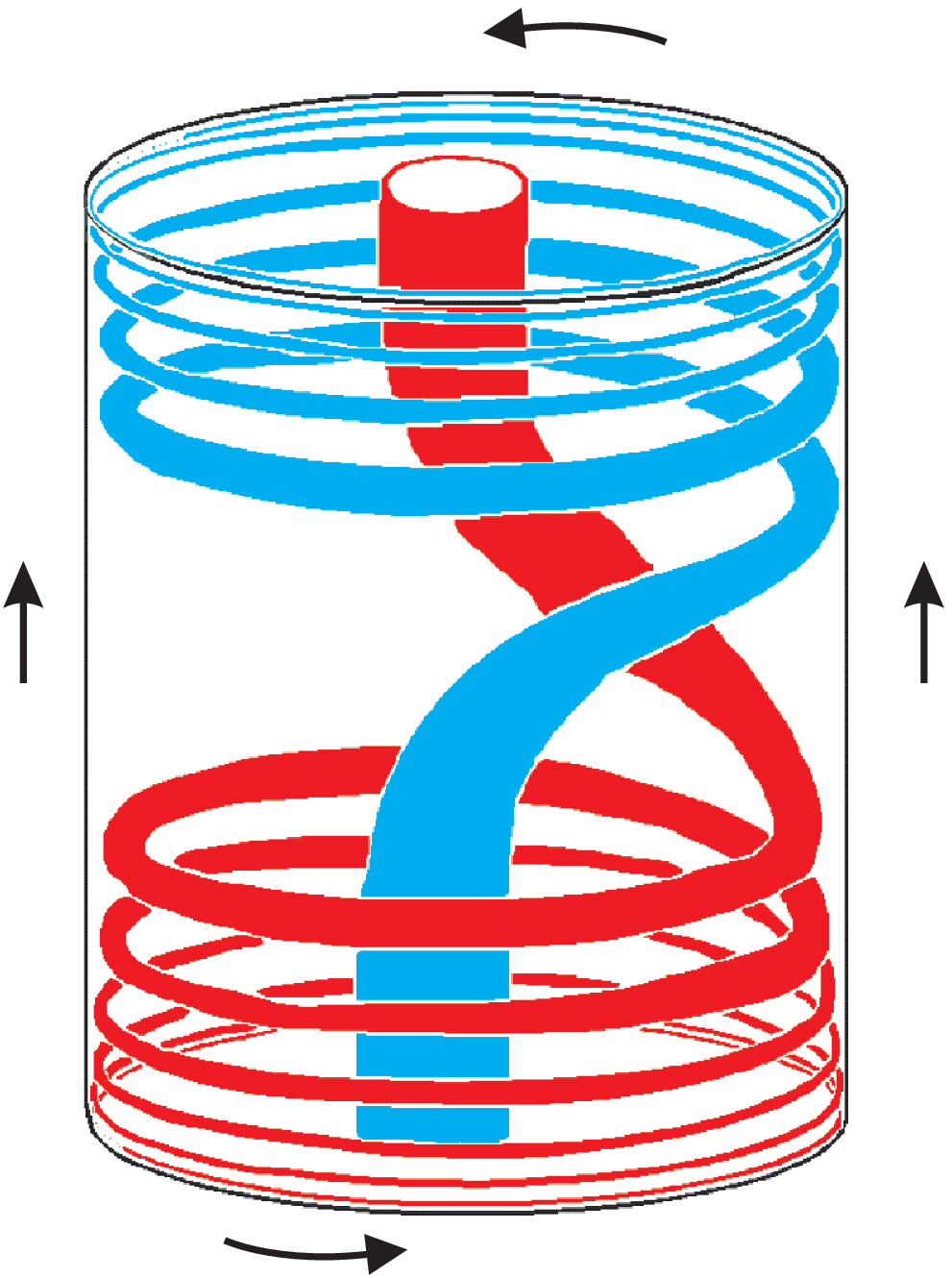}}

  \caption {Borsuk's can $B$, \cite{borsuk}} (see \cite{kuperbergs-minc-reed} for an earlier version of this illustration)
  \label{fig:borsukcan}
\end{figure}

 It is very probable that the roots of the plane fixed point problem go to the famous Brouwer's  fixed-point theorem which states that every $n$ dimensional cube $I^n$ has the fixed point property. In general, the theorem may not hold if $I^n$ is replaced by an $I^n$-like continuum. The first example of an $I^n$-like continuum  without the fixed point property was constructed by Borsuk \cite{borsuk} in 1935. The example consisting of a solid cylinder in $\mathbb{R}^3$ with two spiraling tunnels cut off is illustrated in Fig \ref{fig:borsukcan}. Rotating the example in itself as indicated in the figure forces each vertical level of the cylinder upward except for the bottom and the top disks. So, only the centers of these two disks could stay fixed. But the centers of the bottom and top belong to the two tunnels cut off from the continuum. This ingenious example settled the question of extending  Brouwer's  fixed-point theorem to $I^n$-like continua in $\mathbb{R}^n$ for all $n\ne 2$.
The question for $n=2$ is equivalent to the plane fixed point problem since a plane continuum is $I^2$-like if and only if it does not separate the plane.

The plane fixed point problem has been called the most interesting outstanding problem in plane topology, see for example Bing \cite[p.123]{14} and Hagopian \cite[p. 266]{3}.
Articles \cite{bellamy95, 1,2,3,4,5,6} introduce the problem and give details of progress that has been made since 1930.

The plane fixed point problem has been listed as the first problem among the classical problems in continuum theory, see  Cook, Ingram and Lelek \cite{cook-ingram-lelek}, and Prajs and Charatonik \cite{5}.
No question has attracted more interest from continua theorists than this problem (see Rogers \cite[p. 305]{6}),
and many papers connected the topic have appeared since 1930, including \cite{7,8,9,10,11,12,13,14,15,blokh_etal, borsuk, borsuk54, cartwright-littlewood, fearnley-wright,feuerbacher,16,logan,17,19,20,18,21,24,22,23,hagopian-marsh15,hagopian-marsh-prajs12,25,26,27,28,29,30,manka76,31,32,33,34,35,36,minc-noinvariant-indec,37,38,39,40,41}.
In these papers there were many attempts to solve the problem; some of the authors proved the positive answer under different additional assumptions and some of them tried to construct a counterexample.  In spite of such great effort and many difficult partial results, the fixed point problem is still unsolved, even in the case of planar tree-like continua.

Tree-like continua are very important in the plane fixed point problem. A plane continuum is tree-like if and only if it is 1-dimensional and does not separate the plane. The first example of a tree-like continuum without the fixed point property was constructed by Bellamy \cite{12} in 1980 answering a question by Bing \cite[Question 2, p 122]{14}. This break-through construction was followed by other examples of tree-like continua without the fixed point property \cite{38,39,34,fearnley-wright,35,36,37,hagopian-marsh-prajs12,hagopian-marsh15,logan}.
Applying to those continua and fixed-point-free maps a technique by Fugate and Mohler \cite{16}, we get more tree-like continua without the fixed point property. Some of the examples are atriodic. Even though all these continua appear not to be not planar, it would be very difficult to actually prove this assertion in all atriodic cases. There are simply not many methods of showing that an atriodic tree-like is not planar. In any case, the question whether every plane tree-like continuum has the fixed point property is open, and solving it seems to be crucial to solving the plane fixed point problem.

We now return to Borsuk's can $B\subset \mathbb{R}^3$ to illustrate a plane notion of outchannel.
Let $h$ denote a natural extension of the fixed point map indicated in Figure \ref{fig:borsukcan} to a homeomorphism of $\mathbb{R}^3$ onto itself.
Observe that the red tunnel (the one with the opening at the top of the cylinder) may be called the ``out-tunnel'' since $h$ moves its points away from $B$. Similarly, the other tunnel (the blue one with the opening at the top) is the ``in-tunnel'' since $h$ moves its points towards $B$.
Similar notions of an outchannel and inchannel can be defined in the plane (for a fixed point free map $f$ of a nonseparating plane continuum $X$ which carries the boundary of $X$ minimally into itself).
In 1967-1970, Bell \cite{10}, Sieklucki \cite{40}, and Illidais \cite{28}
independently proved that for every fixed-point-free map  $f$ of a nonseparating plane continuum $X$ into itself, there exists a minimal indecomposable continuum $C\subset\operatorname{bd}(X)$ such that $f(C)=C$. It follows from their proofs that $X$ has an outchannel converging to $C$; see also Brechner and Mayer \cite{brechner-mayer81}.
 Subsequently, Bell noticed that $X$ has exactly one outchanel; for a proof and more results on the subject see \cite{blokh_etal}.

In 1983, Mayer \cite{mayer83} constructed a tree-like planar continuum which might allow for the existence of an outchannel under each embedding in the plane.
This continuum is tree-like, atriodic and indecomposable. It is not weakly chainable, and therefore not chainable (chainable = arc-like).  It is still unknown whether Mayer's continuum has the fixed point property, see \cite[Problem 12 p. 206]{3}. The continuum is a very good illustration of the general weakness of the existing methods in the theory of the plane fixed point problem. The properties of Mayer's continuum are carefully chosen so none of the existing positive partial results can be applied. On the other hand, it is possible that the continuum (or, more likely, its variation) does not have the fixed point property. In that case, how to construct a mapping without a fixed point? Or, better yet, how to systematically approach a construction of a plane tree-like continuum together with a fixed-point-free map?

In 1963 Mioduszewski \cite{mioduszewski} gave necessary and sufficient conditions for a compact metric space to be a continuous image of another one. His conditions are expressed in terms of almost commutative inverse diagrams. In 1994 Feuerbacher \cite[Theorem 3]{feuerbacher} extended Mioduszewski's results to an inverse limit characterization of the fixed-point property for a map from  an arbitrary compact metric space $H=\varprojlim\ofc{X_i,f_i^{i+1}}$ into itself (where $x_i$'s are polyhedra). Mioduszewski-Feuerbacher characterization is general and does not say if $H$ can be embedded in the plane. Possibly, a construction of a nonseparating plane continuum could be accomplished by requiring additionally that $\varprojlim\ofc{X_i,f_i^{i+1}}$ satisfies Anderson-Choquet embedding theorem  \cite[Th. 2.10 p. 23]{4}. However, this could pile an extra condition onto already prohibitively difficult problem since it is not known whether each planar continuum can be re-embedded in $\mathbb{R}^2$ using Anderson-Choquet embedding theorem.

In our paper, we use sequences of finite families of open sets that follow certain patterns to characterize the continua that do not have the fixed point property. 
Sequences of finite families of open sets have already been used in the past (also in terms of open covers) to characterize various topological properties of continua such as chainability, homogeneity, indecomposability, dimension,  and many others.
For example, Bing characterized the pseudo-arc using sequences of finite chains in the plane that follow crooked patterns, for details see \cite{13,bingP,bingPP}.
Also, Hamilton used open covers in \cite{26} to show that every arc-like continuum has the fixed point property. However,  we have found no evidence that a characterization of the continua that do not have the fixed point property has already been obtained in terms of sequences of finite families of open sets.

In section \ref{s3} of the present paper, we give two new characterizations of when a $\mathcal P$-like continuum $X$ does not have the fixed-point property; see Theorem \ref{main}. The characterizations are stated in terms of sequences of open covers of $X$ that follow fixed-point-free patterns; see section \ref{s3}. On the face of it, Theorem \ref{main} does not have anything to do with the plane. However, this is not the case. In each of Corollaries \ref{cor} and \ref{cor2} we consider a sequence $\of{\mathcal{U}_n}_{n=1}^\infty$ of finite collections of open sets in the plane. The sequence satisfies certain recursive conditions guaranteeing that $ \bigcap_{n=1}^\infty\of{\bigcup\mathcal{U}_n}\subset\mathbb{R}^2$ is a tree-like continuum without the fixed point property. On the other hand, if there is no such sequence then every planar tree-like continuum has the fixed point property. Since the conditions imposed on $\of{\mathcal{U}_n}_{n=1}^\infty$ are recursive (the conditions on $\mathcal{U}_n$ depend only on $\mathcal{U}_1,\dots,\mathcal{U}_{n-1}$)
a construction of such a sequence, if it exists, seems quite plausible. This (cautiously) optimistic outlook is further encouraged by our construction, for an arbitrary positive integer $l$, a finite sequence $\of{\mathcal{U}_n}_{n=1}^l$ following the inclusion patterns required in Corollaries \ref{cor} and \ref{cor2}; see Theorem \ref{problem 1} for the statement, the actual construction is given is sections \ref{s4}-\ref{s6}.

\section{Definitions and notation}\label{s2}
 A {\em fixed point\/} of a function $f:X\rightarrow X$ is a point $p$ in $X$ such that $f(p)=p$. A topological space $X$ is said to have the {\em fixed point property\/} if every continuous function from $X$ into $X$ has a fixed point.

 A {\em continuum\/} is a nonempty compact and connected metric space.
A {\em graph\/} is a continuum that can be written as the union of finitely many arcs, any two of which are either disjoint or intersect in one or both of their end points. A graph is {\em acyclic\/} if it does not contain a simple closed curve. An acyclic graph is called a {\em tree\/}.

Let $(X,d)$ be a metric space.  If $r$ is a positive real number, $A$, $A_1$ and $A_2$ are subsets of $X$ and $z\in X$ then we define
 $$
d\of{A,z}=\inf\ofc{d\of{a,z}\mid a\in A},
$$
$$
d\of{A_1,A_2}=\inf\ofc{d\of{a_1,a_2}\mid a_1\in A_1, a_2\in A_2},
$$
 and $B\of{A,r}=\ofc{z\in X\mid d\of{A,z}<r}$.

Let $\mathcal P$ be a class of polyhedra. A finite family $\mathcal{U}$ of open sets in a metric space $X$ is called a $\mathcal P$-{\em cover\/} for $X$ if $\bigcup_{U\in\mathcal U}U=X$ and the geometric realization of the nerve of $\mathcal{U}$ is homeomorphic to a member of $\mathcal P$. The elements of a $\mathcal P$-cover are called {\em links\/}. If $\mathcal P$ is the family of all trees, then a $\mathcal P$-cover is called a {\em tree-cover\/} for $X$.

A family $\mathcal{U}$ of open subsets of a metric space $X$ is a {\em tree-chain\/} if the geometric realization of the nerve of $\mathcal{U}$ is a tree.

A family $\mathcal{U}$ of subsets of a metric space $X$ is {\em taut\/} if for any $U,V\in \mathcal{U}$, $\operatorname{cl}(U)\cap\operatorname{cl}(V)=\emptyset$, if $U\cap V= \emptyset$.

\begin{proposition}\label{p:enlarge}
  Suppose $l$ is a nonnegative integer, $X$ is a compact subset of a metric space $M$ and, for each $n=0,\dots,l$, $\mathcal{U}_n$ is a finite collection of subsets of $X$ such that the collection  $\mathcal{U}=\bigcup_{n=0}^l\mathcal{U}_n$ is taut. Then each $U\in\mathcal{U}$ can be enlarged to a set $U^*$ open in $M$ in such a way that
\begin{enumerate}
  \item $\operatorname{cl}(U^*)\cap\operatorname{cl}(V^*)=\emptyset$ for all $U,V\in \mathcal{U}$ such that $U\cap V=\emptyset$, and
  \item $\operatorname{cl}(U^*)\subseteq V^*$ for all $j=1,\dots,l$, $n=0,\dots,j-1$, $U\in\mathcal{U}_j$, $V\in\mathcal{U}_n$ such that $U\subseteq V$.
\end{enumerate}
\end{proposition}
\begin{proof}
  Let $d$ denote the metric on $M$.  Set 
	$$
	m=\frac{1}{3}\min\ofc{d\of{U,V}\mid U,V\in\mathcal{U},U\cap V=\emptyset}.
	$$
	Finally, set $U^*=B\of{U,2^{-n}m}$ for each $n$ and each $U\in\mathcal{U}_n$, and observe that so defined operation \ ${}^*$ \ satisfies the conclusion of the proposition.
\end{proof}

We say that the {\em mesh} of a family $\mathcal U$ of subsets of a metric space $X$ is the least upper bound of the set
$\{\textup{diam}(U)\ | \ U\in\mathcal U\}$.

%

Let $\mathcal{U}$ and $\mathcal{V}$ be any families of subsets of a metric space $X$. We say that $\mathcal{V}$
\begin{enumerate}
\item {\em refines\/} $\mathcal{U}$ if for any $V\in \mathcal{V}$ there is a $U\in \mathcal{U}$ such that $V\subseteq U$;
\item {\em strongly refines\/}  $\mathcal{U}$ if for any $V\in \mathcal{V}$ there is a $U\in \mathcal{U}$ such that $\operatorname{cl}(V)\subseteq U$.
\end{enumerate}



We also use the following well-known results.
\begin{proposition}\label{lebeg}
Let $X$ be a compact metric space and let $\mathcal U$ and $\mathcal V$ be open covers of $X$ such that the mesh of $\mathcal V$ is less than a Lebesgue number\footnote{A positive number $\lambda$ is a Lebesgue number of an open cover $\mathcal U$ of $X$ if for each set $A\subseteq X$ with diameter less than $\lambda$, there is $U$ in $\mathcal U$ such that $A\subseteq U$.} of $\mathcal U$. Then $\mathcal V$ strongly refines $\mathcal U$.
\end{proposition}
\begin{proposition}\label{taut}
Let $X$ be a continuum. Then for each finite open cover $\mathcal U=\{U_1,U_2,U_3,\ldots ,U_k\}$ of $X$ there is a taut cover $\mathcal V=\{V_1,V_2,V_3,\ldots ,V_k\}$ of $X$ such that $\operatorname{cl}(V_i)\subseteq U_i$ for each $i\in\{1,2,3,\ldots ,k\}$ and the nerves of $\mathcal U$ and $\mathcal V$ are isomorphic.
\end{proposition}
Let $X$ be a continuum. Then $X$ is $\mathcal P$-like if and only if each open cover of $X$ can be refined by a finite open $\mathcal P$-cover of $X$. For details, see \cite[Theorem 5]{mccord}. It follows easily that the following statements are equivalent.
\begin{enumerate}
\item $X$ is $\mathcal P$-like.
\item There is a sequence $(\mathcal{V}_n)$ of taut $\mathcal P$-covers in $X$  such that
\begin{enumerate}
\item for each positive integer $n$, $\mathcal{V}_{n+1}$ strongly refines $\mathcal{V}_n$,
\item for each positive integer $n$, $\mathcal{V}_{n}$ is an $\frac{1}{2^n}$-$\mathcal P$-cover (meaning that each member of $\mathcal{V}_{n}$ has diameter less than $\frac{1}{2^n}$).
\end{enumerate}
\end{enumerate}

%
%
%
%
%


\section{Fixed-point-free patterns of open covers}\label{s3}

Theorem \ref{main} below is the main result in this section. It gives a characterization of $\mathcal P$-like continua that do not have the fixed point property in terms of open covers following special patterns, defined in the following definition.

\begin{definition}\label{pattern}
Let  $(\mathcal{U}_n)$ be a sequence of finite families of subsets of a metric space $X$. Given a sequence $(\varphi_n)$ of functions $\varphi_n:\mathcal{U}_{n+1}\rightarrow \mathcal{U}_n$, consider the following conditions:
\begin{enumerate}
\item[(D1)]\label{D1} for each $U\in \mathcal U_1$ and for each $V\in \mathcal U_2$,
$$
U\cap V\neq \emptyset ~~~ \Longrightarrow ~~~ \varphi_1(V)\cap U=\emptyset,
$$
\item[(D2)]\label{D2} for each positive integer $n$ and for all $U\in \mathcal U_{n+2}$ and $V\in \mathcal U_{n+1}$,
$$
U\cap V\neq \emptyset ~~~ \Longrightarrow  ~~~ \varphi_{n+1}(U)\cap \varphi_n(V)\neq \emptyset,
$$
\item[(D2')]\label{D2'} for each positive integer $n$ and for all $U\in \mathcal U_{n+2}$ and $V\in \mathcal U_{n+1}$,
$$
\operatorname{cl}(U)\subseteq V ~~~ \Longrightarrow  ~~~ \varphi_{n+1}(U)\subseteq \varphi_n(V),
$$
\item[(D3)]\label{D3} for each positive integer $n$ and for all $U,V\in \mathcal U_{n+1}$,
$$
U\cap V\neq \emptyset  ~~~ \Longrightarrow  ~~~ \varphi_n(U)\cap \varphi_n(V)\neq \emptyset.
$$
\end{enumerate}
We say that the sequence $(\mathcal{U}_n)$ of finite families of subsets of a metric space $X$
\begin{itemize}
\item {\em follows a  fixed-point-free pattern}, if there is a sequence $(\varphi_n)$ of functions $\varphi_n:\mathcal{U}_{n+1}\rightarrow \mathcal{U}_n$ satisfying (D1), (D2) and (D3).
\item {\em follows a type-(D2') fixed-point-free pattern}, if there is a sequence $(\varphi_n)$ of functions $\varphi_n:\mathcal{U}_{n+1}\rightarrow \mathcal{U}_n$ satisfying (D1), (D2') and (D3).
\item {\em follows a strong fixed-point-free pattern}, if there is a sequence $(\varphi_n)$ of functions $\varphi_n:\mathcal{U}_{n+1}\rightarrow \mathcal{U}_n$ satisfying (D1), (D2), (D2') and (D3).
\end{itemize}
\end{definition}
\begin{theorem}\label{main}
Let $X$ be a continuum and let $\mathcal P$ be a class of polyhedra. The following statements are equivalent.
\begin{enumerate}
\item[(T1)] \label{(1)} $X$ is a $\mathcal P$-like continuum that does not have the fixed-point property.
\item[(T2)] \label{(2)} There is a sequence $(\mathcal{U}_n)$ of finite open and taut $\mathcal P$-covers of $X$ such that
\begin{enumerate}
\item for each positive integer $n$, $\mathcal{U}_{n+1}$ strongly refines $\mathcal{U}_n$,
\item for each positive integer $n\geq 2$, $\mathcal{U}_{n}$ is a $\frac{\rho}{2^n}$-$\mathcal P$-cover, where
$$
\rho=\min\{d(U,V) \  \vert \ U\cap V=\emptyset, U,V\in \mathcal U_1\}, and
$$
\item $(\mathcal{U}_n)$ follows a fixed-point-free pattern.
\end{enumerate}
\item[(T3)] \label{(3)} There is a sequence $(\mathcal{U}_n)$ of finite open and taut $\mathcal P$-covers of $X$ such that
\begin{enumerate}
\item for each positive integer $n$, $\mathcal{U}_{n+1}$ strongly refines $\mathcal{U}_n$,
\item $\displaystyle \lim_{n\to\infty}\textup{mesh}(\mathcal U_n)=0$, and
\item $(\mathcal{U}_n)$ follows a type-(D2') fixed-point-free pattern.
\end{enumerate}
\end{enumerate}
\end{theorem}

\begin{proof}
We first prove that (T1) implies (T2) and (T3).
Suppose that $f:X\rightarrow X$ is a fixed-point-free map. Since $X$ is compact, there is a positive number $R$ such that $d(x,f(x))>R$ for each $x\in X$.
Let $\mathcal U_1$ be a finite open and taut $\mathcal P$-cover of $X$ with mesh less than $\frac{R}{2}$ such that there are $U,V\in \mathcal U_1$ such that $U\cap V=\emptyset$. Set
$$
\rho=\min\{d(U,V) \  \vert \ U\cap V=\emptyset, U,V\in \mathcal U_1\}.
$$
We construct $\mathcal U_1, \mathcal U_2, \mathcal U_3, \ldots$ by induction. Suppose that $\mathcal U_n$ has been constructed for some integer $n\geq 1$. We now construct $\mathcal U_{n+1}$.

Let $\lambda_n$ be a Lebesgue number for $\mathcal U_n$. Since $f$ is continuous and $X$ is compact, there is $\delta_n>0$ such that $\diam(f(A))<\frac{1}{3}\lambda_n$ for each $A\subseteq X$ with $\diam(A)<\delta_n$. Set
$$
\varepsilon_n=\min\{\frac{1}{3}\lambda_n,\delta_n,\frac{\rho}{2^{n+1}}\}.
$$
We complete our construction of $\mathcal U_1, \mathcal U_2, \mathcal U_3, \ldots$ by letting $\mathcal U_{n+1}$ be a finite open and taut $\mathcal P$-cover of $X$ with mesh less than $\varepsilon_n$.

For each $A\subseteq X$ and $r>0$, let $B(A,r)=\{x\in X \ | \ d(x,A)\leq r\}$. Since $\varepsilon_n\leq \delta_n$ and $\varepsilon_n\leq \frac{1}{3}\lambda_n$,
$$
\diam(B(f(\Cl(V))),\textup{mesh}(\mathcal U_{n+1}))<\lambda_n
$$
for each $V\in \mathcal U_{n+1}$. So, for each $V\in\mathcal U_{n+1}$ there is $\varphi_n(V)\in \mathcal U_n$ such that
$$
	(P_n) \  \  \   \   \  \ \  \  \   \   \  \  B(f(\Cl(V)),\textup{mesh}(\mathcal U_{n+1}))\subseteq \varphi_n(V).
$$

By Proposition \ref{lebeg} $\mathcal U_{n+1}$ strongly refines $\mathcal U_n$ since $\varepsilon_n<\lambda_n$. So (T2)(a) and (T3)(a) are satisfied. (T2)(b) and (T3)(b) follow from the inequality $\varepsilon_n\leq \frac{\rho}{2^{n+1}}$. We now show that $(\mathcal U_n)$ follows a strong fixed-point-free pattern.
\begin{enumerate}
\item[(D1):] Take $U\in \mathcal U_1$ and $V\in \mathcal U_2$ such that $U\cap V\neq \emptyset$. Suppose that $\varphi_1(V)\cap U\neq \emptyset$. Then $\diam(\varphi_1(V)\cup U)<R$. Take $p\in U\cap V$. Then $f(p)\in \varphi_1(V)$ by $(P_1)$. So $d(p,f(p))<R$, a contradiction.
\item[(D2):] Take $U\in \mathcal U_{n+2}$ and $V\in \mathcal U_{n+1}$ such that $p\in U\cap V$ for some $p\in X$.   Then $f(p)\in f(U)\subseteq \varphi_{n+1}(U)$ by $(P_{n+1})$ and  $f(p)\in f(V)\subseteq \varphi_{n}(V)$ by $(P_{n})$. So $\varphi_{n+1}(U)\cap \varphi_n(V)\neq \emptyset$.
\item[(D2'):] Take $U\in \mathcal U_{n+2}$ and $V\in \mathcal U_{n+1}$ such that $U\subseteq V$. Then $f(U)\subseteq \varphi_{n+1}(U)$ by $(P_{n+1})$. It follows that
$$
\varphi_{n+1}(U)\subseteq B(f(\Cl(U)),\textup{mesh}(\mathcal U_{n+1}))\subseteq B(f(\Cl(V)),\textup{mesh}(\mathcal U_{n+1}))\subseteq \varphi_n(V),
$$
 where the last inclusion is $(P_{n})$.
\item[(D3):] Take $U,V\in \mathcal U_{n+1}$ such that $p\in U\cap V$ for some $p\in X$. Then $f(p)\in f(U)\subseteq \varphi_n(U)$ and $f(p)\in f(V)\subseteq \varphi_n(V)$. Therefore, $\varphi_n(U)\cap \varphi_n(V)\neq \emptyset$.
\end{enumerate}

	We have just proved that (T2) and (T3) follow from (T1). We now prove that each of (T2) and (T3) implies (T1). Observe that $X$ is $\mathcal P$-like if either (T2) or (T3) is assumed. So, to complete the proof it is enough to define two fixed-point-free maps $g,h:X\rightarrow X$, one of them assuming (T2) and the other assuming (T3).

	Suppose (T2). A sequence $S=(S_n)$ is called a standard sequence of neighborhoods of $x\in X$ if $x\in S_n\in \mathcal U_n$ for each positive integer of $n$. Notice that there may be many standard sequences of neighborhoods of the same point. For each $x\in X$, we fix one such standard sequence $S^x=(S_n^x)$ of neighborhoods of $x$.
		
\-

{\em Claim.} The following statements are true for all $x\in X$.
\begin{enumerate}
\item $\varphi_{n-1}(S_n^x)\cap \varphi_n(S_{n+1}^x)\neq \emptyset$ for all $n\geq 2$.
\item $\displaystyle \diam\left(\bigcup_{n=i}^{j}\varphi_{n-1}(S_n^x)\right)\leq \sum_{n=i}^{j}\diam(\varphi_{n-1}(S_n^x))$ for all $i$ and $j$ such that $2\leq i\leq j$.
\item $\displaystyle \diam\left(\Cl\left (\bigcup_{n=i}^{\infty}\varphi_{n-1}(S_n^x)\right)\right) < \sum_{n=i}^{\infty}\frac{\rho}{2^{n-1}}=\frac{\rho}{2^{i-2}}$ for all $i\geq 3$.
\item $\displaystyle \bigcap_{i=2}^{\infty}\left(\bigcup_{n=i}^{\infty}\varphi_{n-1}(S_n^x)\right)$ consists of a single point. We denote this point by $g(x)$.
\end{enumerate}		
	
Proof of Claim. (1) follows from (D2). (2) follows by induction from (1) and the triangle inequality. (3) follows from (2) and (T2)(b). (3) implies (4) since $X$ is compact. So, the claim is true.

	\-
	
	We now prove that the function $g:X\rightarrow X$ defined in (4) is continuous. For that purpose consider any $x\in X$ and any $\varepsilon >0$. There is an integer $i\geq 3$ such that $\frac{\rho}{2^{i-2}}<\frac{\varepsilon}{2}$. Take arbitrary $z\in S_i^x$. Since $z\in S_i^x\cap S_i^z$, we infer from $(D3)$ that there exists $p\in \varphi_{i-1}(S_i^x)\cap \varphi_{i-1}(S_i^z)$. It follows from (3) that $d(p,g(x))<\frac{\rho}{2^{i-2}}<\frac{\varepsilon}{2}$ and $d(p,g(z))<\frac{\rho}{2^{i-2}}<\frac{\varepsilon}{2}$. Consequently, $d(g(x),g(z))<\varepsilon$ and $g$ is continuous at each $x\in X$. (Observe that since $g$ is continuous at $x$ for any choice of a standard sequence $S^x=(S_n^x)_{n=2}^{\infty}$ of neighborhoods of $x$, the value of $g(x)$ does not depend on that choice.)
	
	To complete the proof of the implication from (T2) to (T1), it is enough to show that $g(x)\neq x$ for all $x\in X$. It follows from (D2) that there is a point $p\in \varphi_{1}(S_2^x)\cap \varphi_2(S_3^x)$. Since both $p$ and $g(x)$ belong to $\Cl(\bigcup_{n=3}^{\infty}\varphi_{n-1}(S_n^x))$, (3) implies that $d(p,g(x))<\frac{\rho}{2}$. Since $\mathcal U_2$ refines $\mathcal U_1$, $S_2^x$ is contained in some $U\in \mathcal U_1$. By (D1), $\varphi_1(S_2^x)\cap U=\emptyset$. Since both $\varphi_1(S_2^x)$ and $U$ are in $\mathcal U_1$, $p\in \varphi_1(S_2^x)$ and $x\in S_2^x\subseteq U$, it follows from the definition of $\rho$ in (T2)(b) that $d(p,x)\geq \rho$. So, $d(p,g(x))<\frac{\rho}{2}$ implies $g(x)\neq x$.
	
	\-
	
	Suppose (T3). Our construction of a fixed-point-free map $h$ is similar to the construction of $g$, but we need to define standard sequences of neighborhoods differently. A sequence $T=(T_n)$ is called the standard type-(D2') sequence of neighborhoods of $x\in X$ if $x\in T_n\in \mathcal U_n$ and $\Cl(T_{n+1})\subseteq T_n$ for each positive integer $n$.
	
	We will follow here the classic inductive proof of D. K{\"o}nig's infinity lemma to show that each point $x\in X$ has a standard type-(D2') sequence of neighborhoods\footnote{The same result may be also obtained from \cite[Theorem 1]{Cowen}}. Set $N_0=\{k\in\mathbb N \ | \ k\geq 2\}$. For each $k\in N_0$ there is $V_k^{(k)}\in \mathcal U_k$ such that $x\in V_k^{(k)}$. Using repeatedly (T3)(a) we get sets $V_{k-1}^{(k)}$, $V_{k-2}^{(k)}$, $V_{k-3}^{(k)}$, $\ldots$, $V_{1}^{(k)}$ such that $\Cl(V_{n+1}^{(k)})\subseteq V_{n}^{(k)}\in \mathcal U_n$. for each $n=k-1,k-2,k-3,\ldots,1$. We will now construct by induction two sequences $(N_n)$ and $(T_n)$ such that
	$N_n\subseteq \mathbb N$ is an infinite subset of $N_{n-1}$ and
	$T_n\in \mathcal U_n$ is such that $T_n=V_{n}^{(k)}$ for all $k\in N_n$. Suppose that $N_{n-1}$ has already been constructed. Since $N_{n-1}$ is infinite, $\mathcal U_n$ is finite and $V_{n}^{(k)}\in \mathcal U_n$ for all $k$, there is $T_n\in \mathcal U_n$ such that $T_n=V_n^{(k)}$ for infinitely many $k\in N_{n-1}$. To complete the construction of the two sequences, set $N_n=\{k\in N_{n-1} \ | \ T_n=V_n^{(k)}\}$. It is easy to observe that so constructed sequence is a standard type-(D2') sequence of neighborhoods of $x$.
	
	Notice that there may be many standard type-(D2') sequences of neighborhoods of the same point. For each $x\in X$, we fix one such sequence $T^x=(T_n^x)$.
	
	Since $\Cl(T_{n+1}^{x})\subseteq T_n^{x}$, (D2') implies $\varphi_n(T_{n+1}^x)\subseteq \varphi_{n-1}(T_n^x)$ for each $x\in X$. Consequently, $\Cl(\varphi_n(T_{n+1}^x))\subseteq \Cl(\varphi_{n-1}(T_n^x))$ for each $x\in X$ and all $n\geq 2$. It follows from (T3)(b) that $\bigcap_{n=2}^{\infty}\Cl(\varphi_{n-1}(T_n^x))$ consists of a single point. We denote this point by $h(x)$.
	
	Our proof of continuity of $h:X\rightarrow X$ is similar to that for $g$. Consider any $x\in X$ and any $\varepsilon>0$. There is an integer $i\geq 2$ such that $\textup{mesh}(\mathcal U_{i-1})<\frac{\varepsilon}{2}$; see (T3)(b). Take arbitrary $z\in T_i^x$. Since $z\in T_i^x\cap T_i^z$, we infer from $(D3)$ that there exists $p\in \varphi_{i-1}(T_i^x)\cap \varphi_{i-1}(T_i^z)$. Since both $p$ and $h(x)$ both belong to $\Cl(\varphi_{i-1}(T_i^x))$, $d(p,h(x))\leq \textup{mesh}(\mathcal U_{i-1})<\frac{\varepsilon}{2}$. Similarly, $d(p,h(z))<\frac{\varepsilon}{2}$. Consequently, $d(h(x),h(z))<\varepsilon$ and $h$ is continuous at each $x\in X$. (Observe that $h$ is continuous at $x$ for any choice of a standard type-(D2') sequence $T^x=(T_n^x)_{n=2}^{\infty}$ of neighborhoods of $x$, the value of $h(x)$ does not depend on that choice.)
	
To complete the proof of the theorem, it is enough to show that $h(x)\neq x$ for each $x\in X$. Recall that $x\in T_2^x\subseteq T_1^x$, $T_1^x\in \mathcal U_1$, $T_2^x\in \mathcal U_2$ and $h(x)\in \Cl (\varphi_1(T_2^x))$. Using (D1) with $U=T_1^x$ and $V=T_2^x$ we infer that $\varphi_1(T_2^x)\cap T_1^x=\emptyset$. Thus, $\Cl(\varphi_1(T_2^x))\cap T_1^x=\emptyset$ since $T_1^x$ is open. So, $h(x)\neq x$.	
	\end{proof}

The following corollaries easily follow. They give a new possible approach of how to construct a planar tree-like continuum which does not have the fixed-point property.
\begin{corollary}\label{cor}
Let $(\mathcal{U}_n)$ be a sequence of taut tree-chains in $\mathbb R^2$ such that
\begin{enumerate}
\item for each positive integer $n$, $\mathcal{U}_{n+1}$  strongly refines $\mathcal{U}_n$,
\item for each positive integer $n\geq 2$, $\textup{mesh}(\mathcal U_n)<\frac{\rho}{2^n}$, where
$$
\rho=\min\{d(U,V) \  \vert \ U\cap V=\emptyset, U,V\in \mathcal U_1\},
$$
\item $(\mathcal{U}_n)$ follows a fixed-point-free pattern.
\end{enumerate}
Then $X=\bigcap_{n=1}^{\infty}(\bigcup \mathcal U_n)$ is a tree-like continuum in $\mathbb{R}^2$ that does not have the fixed-point property.

If there is no such sequence, then every planar tree-like continuum has the fixed point property.
\end{corollary}
\begin{proof}
The claim follows directly from Theorem \ref{main} since the families $\mathcal V_n=\{U\cap X \ | \ U\in \mathcal U_n\}$ are tree-covers of $X$ satisfying (T2) from Theorem \ref{main}.
\end{proof}

\begin{corollary}\label{cor2}
Let $(\mathcal{U}_n)$ be a sequence of taut tree-chains in $\mathbb R^2$ such that
\begin{enumerate}
\item for each positive integer $n$, $\mathcal{U}_{n+1}$  strongly refines $\mathcal{U}_n$,
\item $\displaystyle \lim_{n\to\infty }\textup{mesh}(\mathcal{U}_{n})=0$, and
\item $(\mathcal{U}_n)$ follows a type-(D2') fixed-point-free pattern.
\end{enumerate}
Then $X=\bigcap_{n=1}^{\infty}(\bigcup \mathcal U_n)$ is a tree-like continuum in $\mathbb{R}^2$ that does not have the fixed point property.

If there is no such sequence, then every planar tree-like continuum has the fixed point property.
\end{corollary}
\begin{proof}
The claim follows directly from Theorem \ref{main} since the families $\mathcal V_n=\{U\cap X \ | \ U\in \mathcal U_n\}$ are tree-covers of $X$ satisfying (T3) from Theorem \ref{main}.
\end{proof}

Next, we present an example of tree-chains $\mathcal U_1$ and $\mathcal U_2$ in $\mathbb R^2$, and a function $\varphi_1:\mathcal U_2\rightarrow \mathcal U_1$ such that
\begin{enumerate}
\item $\mathcal{U}_{2}$ strongly refines $\mathcal{U}_1$,
\item for each $U\in \mathcal U_1$ and for each $V\in \mathcal U_2$,
$$
U\cap V\neq \emptyset ~~~ \Longrightarrow ~~~ \varphi_1(V)\cap U=\emptyset,
$$
and
\item for all $U,V\in \mathcal U_{2}$,
$$
U\cap V\neq \emptyset  ~~~ \Longrightarrow  ~~~ \varphi_1(U)\cap \varphi_1(V)\neq \emptyset .
$$
\end{enumerate}
\begin{example}\label{ex1}
Let $\mathcal U_1=\{U_1,U_2,U_3,\ldots,U_{13}\}$ and $\mathcal U_2=\{V_1,V_2,V_3,\ldots,V_{131}\}$ be tree-chains in $\mathbb R^2$ as pictured on Figure \ref{fig:tree-chains}.

\begin{figure}[h!]
	\centering
		\includegraphics[width=25.0em]{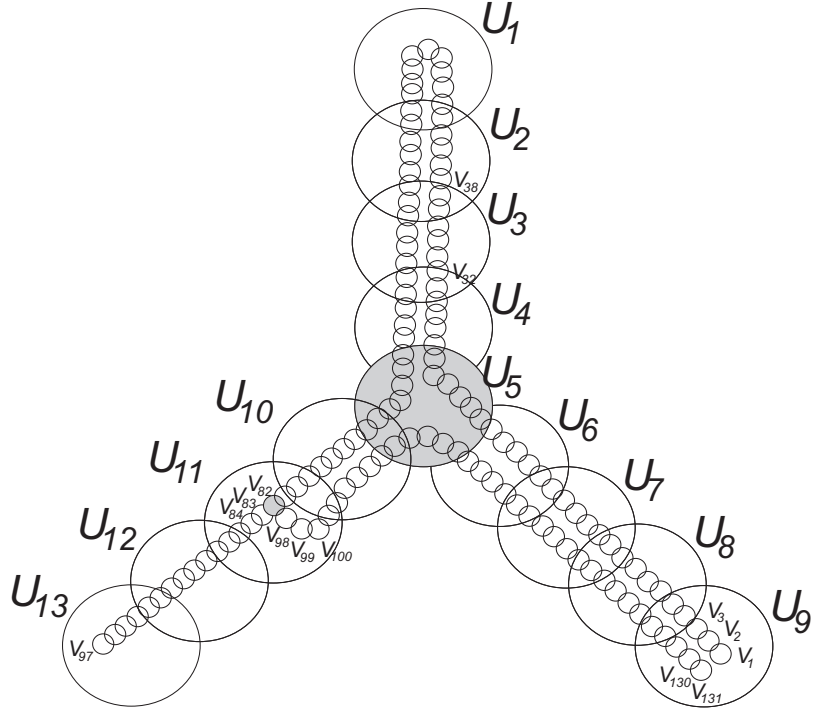}
	\caption{Tree-chains $\mathcal U_1$ and $\mathcal U_2$}
	\label{fig:tree-chains}
\end{figure}

We define $\varphi_1:\mathcal U_2\rightarrow\mathcal U_1$ by
\begin{itemize}
\item $\varphi_1(V_i)=U_{13}$ for each $i\in \{1,2,3,\ldots,32\}$,
\item $\varphi_1(V_{33})=U_{12}$, $\varphi_1(V_{34})=U_{11}$, $\varphi_1(V_{35})=U_{10}$, $\varphi_1(V_{36})=U_{5}$, $\varphi_1(V_{37})=U_{6}$,
\item $\varphi_1(V_{i})=U_{7}$ for each $i\in \{38, 39, 40, \ldots, 82\}$,
\item $\varphi_1(V_{83})=U_{8}$,
\item $\varphi_1(V_{i})=U_{9}$ for each $i\in \{84, 85, 86, \ldots, 97\}$,
\item $\varphi_1(V_{98})=U_{6}$, $\varphi_1(V_{99})=U_{5}$, $\varphi_1(V_{100})=U_{4}$, $\varphi_1(V_{101})=U_{3}$, $\varphi_1(V_{102})=U_{2}$,
\item $\varphi_1(V_{i})=U_{1}$ for each $i\in \{103, 104, 105, \ldots, 131\}$.
\end{itemize}
It is obvious that the tree-chains $\mathcal U_1$ and $\mathcal U_2$, and the function $\varphi_1$ satisfy all the required properties.
\end{example}

\begin{problem}
Does there exist a sequence of tree-chains in $\mathbb R^2$ satisfying all the conditions from Corollary \ref{cor} or Corollary \ref{cor2}?
\end{problem}

A negative answer to the above problem would prove that every planar tree-like continuum has the fixed point property. On the other hand, a positive answer (i.e. constructing an infinite sequence of tree-chains  $\mathcal U_1, \mathcal U_2,\dots$ satisfying the conditions in either \ref{cor} or \ref{cor2}) would yield an example of a tree-like continuum without the fixed point property. In Example \ref{ex1},  we constructed a sequence of just two tree-chains $\mathcal U_1, \mathcal U_2$ satisfying the required conditions. We will extend this construction to get a sequence of arbitrary finite length. More precisely, we will prove the following theorem.

\begin{theorem}\label{problem 1}
For every positive integer $n$, there exist tree-chains $\mathcal U_1$, $\mathcal U_2$, $\ldots $, $\mathcal U_n$ in $\mathbb R^2$, and functions $\varphi_1:\mathcal U_2\rightarrow \mathcal U_1$, $\varphi_2:\mathcal U_3\rightarrow \mathcal U_2$, $\ldots $,  $\varphi_{n-1}:\mathcal U_n\rightarrow \mathcal U_{n-1}$ such that
\begin{enumerate}
\item for each $k\in \{1,2,\ldots ,n-1\}$, $\mathcal{U}_{k+1}$  strongly refines $\mathcal{U}_k$,
\item for each $U\in \mathcal U_1$ and for each $V\in \mathcal U_2$,
$$
U\cap V\neq \emptyset ~~~ \Longrightarrow ~~~ \varphi_1(V)\cap U=\emptyset,
$$
\item for each $k\in \{1,2,\ldots ,n-2\}$ and for all $U\in \mathcal U_{k+2}$ and $V\in \mathcal U_{k+1}$,
$$
U\cap V\neq \emptyset ~~~ \Longrightarrow  ~~~ \varphi_{k+1}(U)\cap \varphi_k(V)\neq \emptyset
$$
or,  for each $k\in \{1,2,\ldots ,n-2\}$ and for all $U\in \mathcal U_{k+2}$ and $V\in \mathcal U_{k+1}$,
$$
\operatorname{cl}(U)\subseteq V ~~~ \Longrightarrow  ~~~ \varphi_{k+1}(U)\subseteq \varphi_k(V),
$$
\item for each $k\in \{1,2,\ldots ,n-1\}$ and for all $U,V\in \mathcal U_{k+1}$,
$$
U\cap V\neq \emptyset  ~~~ \Longrightarrow  ~~~ \varphi_k(U)\cap \varphi_k(V)\neq \emptyset.
$$
\end{enumerate}
\end{theorem}

The main goal of the rest of the paper is to provide a proof to the above theorem. The proof is completed at the end of section \ref{s6}.

\section{Simplicial graphs, simplicial maps, and simplicial diagrams} \label{s4}
In this section we present the concepts of simplicial graphs, maps and diagrams which will be used later to prove Theorem \ref{problem 1}. We also point out several results that are either well-known or easy to prove. We state these results in the form of observations and we leave the easy proofs to the reader.  We begin with the definition of coincident points.

If $f$ and $g$ are two mappings of a space $X$ into a space $Y$, we say that $x\in X$ is a coincidence point of $f$ and $g$ if $f\of{x}=g\of{x}$.

Next we consider a special kind of commutative diagram.  Let $l$ be either a positive integer or $\infty$. Suppose
$$
\xymatrix@1@C=10pt{ X_0 & X_1\ar[l]_{\rule{5pt}{0pt}g_0} & X_2\ar[l]_{\rule{5pt}{0pt}g_1} & \dots\ar[l]_{\rule{5pt}{0pt}g_2} }
$$
is an inverse sequence where $g_n$ is defined for each nonnegative integer $n<l$. Suppose also $i$ and $j$ are integers such that $0\le i<j-1<l$. Then the composition
$g_i\circ g_{i+1}\circ\dots\dots\circ g_{j-1}$ mapping $X_j$ to $X_i$ is denoted by $g_{ij}$. Additionally, let $g_{ii}$ denote the identity on $X_i$.
Sometimes we write $g_{i,j}$ instead of $g_{ij}$, especially when an arithmetic operation is involved. For instance, we write  $g_{i+1,j-1}$ instead of $g_{i+1j-1}$.
When any letter (not necessarily ``g'') is used to denote the bonding maps in a similar context, the same letter with a subscript ''$ij$'' denotes the composition of the bonding maps taking the $j$-th factor space to the $i$-th.

Let $l$ be again either a positive integer or $\infty$. Suppose
$\xymatrix@1@C=10pt{ X_0 & X_1\ar[l]_{\rule{5pt}{0pt}g_0} & X_2\ar[l]_{\rule{5pt}{0pt}g_1} & \dots\ar[l]_{\rule{5pt}{0pt}g_2} }$
and $\xymatrix@1@C=10pt{ X_0 & X_1\ar[l]_{\rule{5pt}{0pt}f_0} & X_2\ar[l]_{\rule{5pt}{0pt}f_1} & \dots\ar[l]_{\rule{5pt}{0pt}f_2} }$
are two inverse sequences where $g_n$ and $f_n$ are defined for each nonnegative integer $n<l$. Then $D_l\of{X_n,g_n,f_n}$ denotes the following diagram:

$$\xymatrixcolsep{2pc}\xymatrix@1{
&X_1\ar[dl]_{f_0}&X_2\ar[dl]_{f_1}\ar[l]_{g_1}&X_3\ar[dl]_{f_2}\ar[l]_{g_2}&X_4\ar[dl]_{f_3}\ar[l]_{g_3}&\dots\ar[l]\\
X_0&X_1\ar[l]_{g_0}&X_2\ar[l]_{g_1}&X_3\ar[l]_{g_2}&X_4\ar[l]_{g_3}&\dots\ar[l]
}
$$

In this paper we are interested in commutative diagrams (meaning that $f_{i-1}\circ g_i=g_{i-1}\circ f_i$ for each positive integer $i$) of continuous surjections with no coincidence points. We say that $D_l\of{X_n,g_n,f_n}$ has no coincidence points if $g_n$ and $f_n$ have no coincidence points for all nonnegative integers $n<l$. We say that a diagram $D_l\of{X_n,g_n,f_n}$ is surjective if all maps $g_n$ are surjections.

\begin{observation}
  Suppose $D_l\of{X_n,g_n,f_n}$ is a commutative diagram such that $f_0$ and $g_0$  have no coincidence points. Then $f_n$ and $g_n$  have no coincidence points for all nonnegative integers $n<l$.
\end{observation}

Next we define simplicial graphs and their geometric realizations.
By a \emph{simplicial graph} $G$ we understand an abstract one-dimensional, finite simplicial complex, that consists of a fixed finite set of vertices $\mathcal{V}\of{G}$ and a fixed set of edges $\mathcal{E}\of{G}$ which consists of a collection of two member sets $\{u,v\}$ where $u,v\in \mathcal{V}\of{G}$ and $u\neq v$. Two vertices belonging to an edge are called \emph{adjacent}. We say that $v,v^{\prime}\in\mathcal{V}\of{G}$ are \emph{$k$-close in $G$} if there are $v_0,v_1,\dots,v_k\in\mathcal{V}\of{G}$ such that $v_0=v$, $v_k=v^{\prime}$, and, for each $n=1,\dots,k$, $v_{n-1}$ and $v_n$ are either equal or adjacent. Notice that $v$ and $v^{\prime}$ are $1$-close if and only if they are either equal or adjacent.

A vertex $v\in \mathcal G$ is an end-point of the graph $G$ if there is only one vertex $u\in \mathcal V(G)$ such that $\{u,v\}\in E(G)$.

If $p,q\in\mathbb{R}^m$, then $\ofi{p,q}$ denotes the straight line segment joining $p$ and $q$ in $\mathbb{R}^n$.

We say that an injection $i:\mathcal{V}\of{G}\to\mathbb{R}^m$ is \emph{consistent with $G$} if the following two conditions are satisfied for all vertices $u$ and $v$  adjacent in $G$:
\begin{itemize}
  \item $\ofi{i(u),i(v)}\cap i\of{\mathcal{V}(G)}=\ofc{i(u),i(v)}$, and
 \item $\ofi{i(u),i(v)}\cap\ofi{i(u^{\prime}),i(v^{\prime})}\subset\ofc{i(u),i(v),i(u^{\prime}),i(v^{\prime})}$ for all edges $\ofc{u^{\prime},v^{\prime}}\in\mathcal{E}(G)$ such that $\ofc{u^{\prime},v^{\prime}}\ne\ofc{u,v}$.
\end{itemize}

\begin{observation}
  For each simplicial graph $G$ there exists an injection $i:\mathcal{V}\of{G}\to\mathbb{R}^3$ which is consistent with $G$.
\end{observation}

Suppose $G$ is a simplicial graph and $i:\mathcal{V}\of{G}\to\mathbb{R}^m$ is an injection consistent with $G$.  Let $\ofa{G}_i$ denote the union of $i\of{\mathcal{V}(G)}$ and all segments $\ofi{i(u),i(v)}$ where $u$ and $v$ are adjacent in $G$. In this context we say that $\ofa{G}_i$ is the \emph{geometric realization} of $G$ associated with $i$.
If $\mathcal{V}\of{G}\subset\mathbb{R}^m$ and $i$ is the inclusion, we will simply write $\ofa{G}$ instead of $\ofa{G}_i$.

 Any two geometric realizations of a simplicial graph $G$ are isomorphic in the following sense.
\begin{observation}
  Suppose $G$ is a simplicial graph. Let $i:\mathcal{V}\of{G}\to\mathbb{R}^m$ and $j:\mathcal{V}\of{G}\to\mathbb{R}^{n}$. Then there is a homeomorphism $h$ of $\ofa{G}_i$ onto $\ofa{G}_j$ such that $h\of{i(v)}=j(v)$ for each $v\in\mathcal{V}\of{G}$, and $h$ is linear on $\ofi{i(u),i(u)}$ for all $u$ and $v$ adjacent in $G$.
\end{observation}

We say that $G$ is connected if its geometric realization $\ofa G$ is connected. Similarly, $G$ is a (simplicial) tree if its geometric realization $\ofa G$ is tree (a connected union of finitely many arcs with no simple closed curve).

To simplify the notation we may assume without loss of generality that the set of vertices of an arbitrary simplicial graph $G$ is contained in $\mathbb{R}^m$ (where $m$ is $3$ or $2$ if possible) in such a way that the inclusion is consistent with $G$.
Then each abstract edge $\ofc{u,v}$ uniquely corresponds to the geometric straight linear edge $\ofi{u,v}\subset\ofa{G}\subset\mathbb{R}^m$.  Observe that if two different geometric segments in $G$ intersect, their intersection is a common vertex of $G$.
Moreover, each component of $\ofa{G}\setminus\mathcal{V}\of{G}$ is the interior of $\ofi{u,v}$ for some $\ofc{u,v}\in\mathcal{E}\of{G}$.



 A \emph{simplicial map} of a simplicial graph $G$ into a simplicial graph $H$ is a function $f:\mathcal{V}\of{G}\to\mathcal{V}\of{H}$ taking every two vertices adjacent in $G$  either onto a pair of  vertices adjacent in $H$ or onto a single vertex.

\begin{observation}\label{o:fkclose}
  Suppose $G$ and $H$ are simplicial graphs and $f:\mathcal V\of G\to \mathcal V\of H$ is a simplicial map. Then  $f\of{v}$ and $f\of{v^{\prime}}$ are $k$-close in $H$ for all positive integers $k$ and all vertices $v$ and $v^{\prime}$ that are $k$-close in $G$.
\end{observation}

 Let $G$ be a simplicial graph and let $\ofa G$ be its geometric realization. Note that for each $x\in \ofa G$ either $x\in \mathcal V(G)$ or $x\not \in \mathcal V(G)$. In latter case there are $\{u,v\}\in \mathcal E(G)$ and $t\in (0,1)$ such that $x=(1-t)u+tv$.

 Let $f:\mathcal{V}\of{G}\to\mathcal{V}\of{H}$ be a simplicial map. The function $\ofa f:\ofa G\rightarrow \ofa H$, defined by
 $$
 \ofa f(x)=\begin{cases}
				f(x)\text{;} & x \in \mathcal{V}\of{G}, \\
				(1-t)f(u)+tf(v)\text{;} & x=(1-t)u+tv \text{ for some } \{u,v\}\in \mathcal E(G) \text{ and } t\in (0,1),
			\end{cases}
 $$
is called the geometric realization of the simplicial map $f$.

 Notice that $\ofa f$ is always a continuous mapping.

We say that $f:\mathcal V\of G\to \mathcal V\of H$ is a \emph{simplicial surjection} if $f$ is a simplicial map and $\mathcal{V}\of{H}=f\of{\mathcal{V}\of{G}}$.

 Note that there is a simplicial surjection $f:\mathcal V\of G\to \mathcal V\of H$ such that its geometric realization $\ofa f:\ofa G\rightarrow \ofa H$ is not surjective. However, if $G$ is connected and $H$ is a tree, then $\ofa f$ is  surjective for any simplicial surjection $f:\mathcal V\of G\to \mathcal V\of H$; see the following observation.

\begin{observation}\label{o:simpsurtree}
  Suppose $C$ is a connected simplicial graph,  $T$ is a simplicial tree, and $f:\mathcal V\of C\to \mathcal V\of T$ is a simplicial surjection. Then  $\mathcal{E}\of{T}=f\of{\mathcal{E}\of{C}}$.
In particular, if $u$ and $v$ are $1$-close vertices of $T$, then there are $1$-close vertices $a,b\in\mathcal{V}\of C$ such that $f\of{a}=u$ and $f\of{b}=v$.
\end{observation}

\begin{observation}\label{ob:simpnocoi}
  Let $G$ and $H$ be simplicial graphs, and let  $f, g:\mathcal V\of G\to \mathcal V\of H$ be two simplicial maps. Then $\ofa f$ and $\ofa g$  have no coincidence points (on the geometric realization $\ofa G$ of $G$) if and only if
  \begin{enumerate}
    \item $f\of{v}\ne g\of{v}$ for all $v\in\mathcal{V}\of{G}$, and
    \item $\ofc{f\of{u},f\of{v}}\not\subset\ofc{g\of{u},g\of{v}}$ for all edges $\{u, v\}\in\mathcal{E}\of{G}$.
  \end{enumerate}
\end{observation}

Let $G$ and $H$ be simplicial graphs, and let $f,g:\mathcal V\of G\to \mathcal V\of H$ be simplicial maps. We say that $v\in\mathcal{V}\of{G}$ is a \emph{proximity vertex} for $f$ and $g$ if $f\of{v}$ and $g\of{v}$ are $2$-close in $G$.

\begin{observation}
  Let $G$ and $H$ be simplicial graphs, and let  $f, g:\mathcal V\of G\to \mathcal V\of H$ be two simplicial maps with no proximity vertices. Then $\ofa f$ and $\ofa g$ have no coincidence points.
\end{observation}

We say that a diagram $D_l\of{\mathcal V\of{G_n},g_n,f_n}$ is simplicial if all $G_n$ are simplicial graphs and all $g_n:\mathcal V\of{G_{n+1}} \rightarrow \mathcal V\of{G_n}$ and $f_n:\mathcal V\of{G_{n+1}} \rightarrow \mathcal V\of{G_n}$ are simplicial maps. Let $D_l\of{\mathcal V\of{G_n},g_n,f_n}$ be a simplicial diagram. We say that $D_l\of{\ofa{G_n},\ofa{g_n},\ofa{f_n}}$ is a geometric realization of $D_l\of{\mathcal V\of{G_n},g_n,f_n}$.

\begin{observation}\label{o:nopoximityforfnandgn}
  Suppose $D_l\of{\mathcal V\of{G_n},g_n,f_n}$ is a simplicial commutative diagram such that $f_0$ and $g_0$  have no proximity vertex. Then $f_n$ and $g_n$  have no proximity vertices for all nonnegative integers $n<l$.
\end{observation}

Let $G$ be an arbitrary simplicial graph in $\mathbb R^3$. Let $G^{\of{3}}$ be the graph resulting from subdividing each edge of $G$ into three congruent parts. More precisely, for each edge $e=\{u,v\}\in \mathcal{E}\of{G}$, let $u_e=\frac23 u+\frac13 v$ and $v_e=\frac13 u+\frac23 v$. Let $\mathcal{V}^3(G)$ denote the set of all points $v_e$ where $v\in\mathcal{V}\of{G}$, $e\in\mathcal{E}\of{G}$ and $v\in e$. Set $\mathcal{V}(G^{\of{3}})=\mathcal{V}\of{G}\cup\mathcal{V}^3(G)$. Finally, let $\mathcal{E}(G^{\of{3}})$ be the set of all edges $\{u,u_e\}$, $\{u_e,v_e\}$ and $\{v_e,v\}$ where $e=\{u,v\}\in \mathcal{E}\of{G}$.
Notice that $\ofa{G^{\of{3}}}=\ofa{G}$ for any simplicial graph $G$.
\begin{observation}
Let $G$ be an arbitrary simplicial graph. Then,  at least one of any two adjacent vertices of $G^{\of{3}}$ must belong to $\mathcal{V}^3(G)$.
\end{observation}

\begin{observation}\label{ob:adjacentinG^3}
Let $G$ be an arbitrary simplicial graph, and let $w\in\mathcal{V}^3(G)$. Then there is a unique vertex $u\in\mathcal{V}\of{G}$ and a unique edge $e=\{u,v\}\in \mathcal E\of G$ such that $w\in e$. Consequently, $w$ is either $u_e$ or $v_e$.
\end{observation}

\begin{observation}\label{o:not2closeinG^3}
Let $G$ be an arbitrary simplicial graph. Let $A$ and $B$ two disjoint subsets of its geometric realization $\ofa G$, each of which is either the segment $\ofi{u,v}$ for some edge $\{u,v\}\in \mathcal E \of G$ or $\{v\}$ for some  vertex $v\in \mathcal V \of G$. Finally, let $a\in A\cap \mathcal{V}(G^{\of{3}})$ and $b\in B\cap \mathcal{V}(G^{\of{3}})$. Then $a$ and $b$ are not $2$-close in $G^{\of{3}}$.
\end{observation}

\begin{observation}\label{o:2closeinG^3}
Let $G$ be an arbitrary simplicial graph and let $t$ be either $1/3$ or $2/3$. Let $\{x,y\},\{y,z\}\in \mathcal E \of G$ be two distinct edges of $G$. Then $(1-t)y+tx$ and $ty+(1-t)z$ both belong to $\mathcal{V}^3(G)$, and are not $2$-close in $G^{\of{3}}$.
\end{observation}

Suppose $f:\mathcal V\of G\rightarrow \mathcal V\of H$ is a simplicial map of a graph $G$ to a graph $H$. Define a function $f^{\of{3}}:\mathcal{V}\of{G^{\of{3}}}\to \mathcal{V}\of{H^{\of{3}}}$ by setting $f^{\of{3}}\of{v}=\ofa f\of{v}$ for each $v\in\mathcal{V}(G^{\of{3}})$. Observe that so defined $f^{\of{3}}$ is a simplicial map of $G^{\of{3}}$ into $H^{\of{3}}$.
Notice that $\ofa f:\ofa{G}\rightarrow \ofa{H}$ and $\ofa{f^{\of{3}}}:\ofa{G^{\of{3}}}\rightarrow \ofa{H^{\of{3}}}$ are identical functions.
\begin{observation}
  Suppose $f,g:\mathcal V\of G\to \mathcal V\of H$ are simplicial maps. Let $p$ be an arbitrary point of $\ofa G$. Then $p$ is a coincidence point of $\ofa{f^{\of{3}}}$ and $\ofa{g^{\of{3}}}$ if and only if it is a coincidence point of $\ofa f$ and $\ofa g$.
\end{observation}

\begin{proposition}
  Suppose $G$ and $H$ are simplicial graphs, and  $f, g:\mathcal V\of G\rightarrow \mathcal V\of H$ are two simplicial maps . If $\ofa f:\ofa G\rightarrow \ofa H$ and $\ofa g:\ofa G\rightarrow \ofa H$ have no coincidence points, then $f^{\of{3}}:G^{\of{3}}\to H^{\of{3}}$ and $g^{\of{3}}:G^{\of{3}}\to H^{\of{3}}$ have no proximity vertex.
\end{proposition}

\begin{proof}
Let $w\in\mathcal{V}\of{G^{\of{3}}}$. We need to show that $f^{\of{3}}\of{w}$ and $g^{\of{3}}\of{w}$ are not $2$-close in $H^{\of 3}$.

Suppose $w\in\mathcal{V}(G)$. Then $f^{\of{3}}\of{w}=f\of{w}$ and $g^{\of{3}}\of{w}=g\of{w}$ are two distinct vertices of $H$, and therefore are not $2$-close in $H^{\of{3}}$ by \ref{o:not2closeinG^3}. So, we may assume that $w\in\mathcal{V}^3(G)$. If follows from \ref{ob:adjacentinG^3} that there is an edge $e=\{u,v\}\in \mathcal E\of G$ such that $w=(1-t)u+tv\in \ofi{u,v}$ where $t$ is either $1/3$ or $2/3$. So, $f^{\of{3}}\of{w}\in \ofa{f^{\of{3}}}\of{\ofi{u,v}}=\ofa f\of{\ofi{u,v}}$ and $g^{\of{3}}\of{w}\in \ofa{g^{\of{3}}}\of{\ofi{u,v}}=\ofa g\of{\ofi{u,v}}$. Observe that each of $\ofa f\of{\ofi{u,v}}$ and $\ofa g\of{\ofi{u,v}}$ is either a segment (either $\ofi{f(u),f(v)}$ or $\ofi{g(u),g(v)}$) in $\ofa H$ or a one vertex subset $\{f(v)\}$ or $\{g(v)\}$ of $\ofa H$.
If  $\ofa f\of{\ofi{u,v}}\cap \ofa g\of{\ofi{u,v}}=\emptyset $ then \ref{o:not2closeinG^3} implies that $f^{\of{3}}\of{w}$ and $g^{\of{3}}\of{w}$ are not $2$-close in $H^{\of{3}}$, and the proposition is true. So, we may assume that $\ofa f\of{\ofi{u,v}}\cap \ofa g\of{\ofi{u,v}}\neq\emptyset$.

It follows from \ref{ob:simpnocoi} that either $f\of{u}=g\of{v}$, $f\of{v}\notin\{ g\of{u},g\of v\}$ and $g\of{u}\notin \{ f\of{u},f\of v\}$, or $f\of{v}=g\of{u}$, $f\of{u}\notin \{ g\of{u},g\of v\}$ and $g\of{v}\notin \{ f\of{u},f\of v\}$. Since the situation is symmetric, we may assume that $f\of{u}=g\of{v}$, $f\of{v}\notin\{ g\of{u},g\of v\}$ and $g\of{u}\notin \{ f\of{u},f\of v\}$.
In particular, $f\of{v}\ne g\of{u}$, $f\of{u}\ne g\of{u}$, $f\of{v}\ne g\of{v}$, $f\of{u}\ne f\of{v}$ and $g\of{u}\ne g\of{v}$.
Hence,  $f^{\of{3}}\of{w}=(1-t)f\of{u}+tf\of{v}\in \ofi{f(u),f(v)}$, $g^{\of{3}}\of{w}=(1-t)g\of{u}+tg\of{v}=tf\of{u}+(1-t)g\of{u}\in \ofi{g(u),g(v)}$ where $\{f(u),f(v)\}$ and $\{g(u),g(v)\}$ are two edges of $H$ that meet only in one vertex $f\of{u}=g\of{v}$.
Thus, $f^{\of{3}}\of{w}$ and $g^{\of{3}}\of{w}$ are not $2$-close in $H^{\of 3}$ by \ref{o:2closeinG^3}.
\end{proof}

\begin{corollary} \label{t:nopoximitydiagram}
    Suppose that $D_l\of{\mathcal V\of{G_n},g_n,f_n}$ is a simplicial, surjective and commutative diagram such that $\ofa{g_0}$ and $\ofa{f_0}$  have no coincidence points. Then $D_l\of{\mathcal V\of{G_n^{\of{3}}},g_n^{\of{3}},f_n^{\of{3}}}$  is a simplicial, surjective and commutative diagram such that $g_n^{\of 3}$ and $f_n^{\of 3}$ have no proximity vertices for all nonnegative integers $n<l$. Also, $\ofa{G_n^{\of{3}}}$ and $\ofa{G_n}$ are homeomorphic for all $n\le l$.
\end{corollary}

\section{Long diagrams of trees with no coincidence points}\label{s5}

It is well-known that each infinite commutative diagram $D_\infty\of{X_n,g_n,f_n}$ induces a map $f:\varprojlim \of{X_n,g_n}\to\varprojlim \of{X_n,g_n}$ defined by $f\of{x_n}=\of{f_{n-1}\of{x_n}}$ for all $\of{x_n}\in\varprojlim \of{X_n,g_n}$. The induced map $f$ has no fixed points if and only if $f_n$ and $g_n$ have no coincidence points for some $n$. In 1982, Oversteegen and Rogers \cite{39} constructed a tree-like continuum with a fixed-points-free map induced from an infinite commutative diagram $D_\infty\of{X_n,g_n,f_n}$ with no coincidence points. Another such example was recently constructed by Hernández-Gutiérrez and  Hoehn \cite{logan}.

 \begin{theorem}[Oversteegen and Rogers, and Hernández-Gutiérrez and  Hoehn]\label{t:infdiagramoftrees}
  There exists a commutative, surjective and simplicial diagram \newline $D_{\infty}\of{X_n,g_n,f_n}$ where $X_n$'s are simplicial trees and $f_n$ and $g_n$ have no coincidence points for each $n$.
\end{theorem}

\begin{figure}\label{fig:Tn}
\centering
\begin{picture}(280,165)
\put(11,10){\scalebox{0.5}{\includegraphics{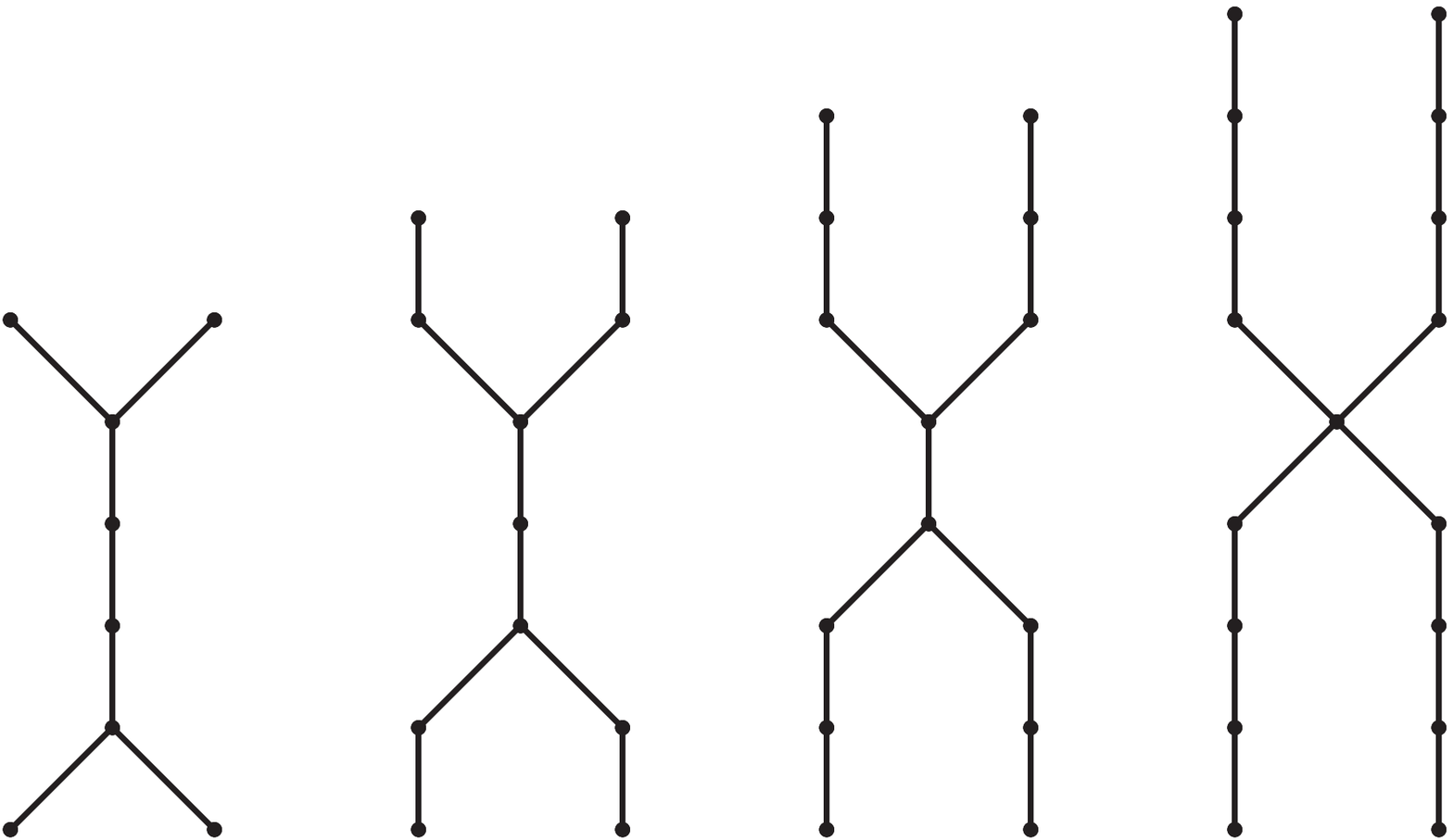}}}
\put(27,0){$T_0^k$}
\put(0,11){$v_{0}^{-1}$}
\put(51,11){$v_{0}^{1}$}
\put(0,101){$v_{5}^{-1}$}
\put(51,101){$v_{5}^{1}$}

\put(34,29){$v_{1}^{0}$}
\put(34,47){$v_{2}^{0}$}
\put(34,65){$v_{3}^{0}$}
\put(34,82){$v_{4}^{0}$}


\put(102,0){$T_1^k$}
\put(88,11){$v_{0}^{-1}$}
\put(88,28){$v_{1}^{-1}$}

\put(124,11){$v_{0}^{1}$}
\put(124,28){$v_{1}^{1}$}

\put(88,101){$v_{5}^{-1}$}
\put(124,101){$v_{5}^{1}$}

\put(88,119){$v_{6}^{-1}$}
\put(124,119){$v_{6}^{1}$}

\put(106,47){$v_{2}^{0}$}
\put(106,65){$v_{3}^{0}$}
\put(106,82){$v_{4}^{0}$}


\put(173,0){$T_2^k$}
\put(161,11){$v_{0}^{-1}$}
\put(161,28){$v_{1}^{-1}$}
\put(161,46){$v_{2}^{-1}$}
\put(196,11){$v_{0}^{1}$}
\put(196,28){$v_{1}^{1}$}
\put(196,46){$v_{2}^{1}$}

\put(161,101){$v_{5}^{-1}$}
\put(196,101){$v_{5}^{1}$}

\put(161,119){$v_{6}^{-1}$}
\put(196,119){$v_{6}^{1}$}

\put(161,137){$v_{7}^{-1}$}
\put(196,137){$v_{7}^{1}$}

\put(179,65){$v_{3}^{0}$}
\put(179,82){$v_{4}^{0}$}


\put(245,0){$T_3^k$}
\put(233,11){$v_{0}^{-1}$}
\put(233,28){$v_{1}^{-1}$}
\put(233,46){$v_{2}^{-1}$}
\put(233,64){$v_{3}^{-1}$}
\put(268,11){$v_{0}^{1}$}
\put(268,28){$v_{1}^{1}$}
\put(268,46){$v_{2}^{1}$}
\put(268,64){$v_{3}^{1}$}

\put(233,101){$v_{5}^{-1}$}
\put(268,101){$v_{5}^{1}$}

\put(233,119){$v_{6}^{-1}$}
\put(268,119){$v_{6}^{1}$}

\put(233,137){$v_{7}^{-1}$}
\put(268,137){$v_{7}^{1}$}

\put(233,155){$v_{8}^{-1}$}
\put(268,155){$v_{8}^{1}$}

\put(251,82){$v_{4}^{0}$}

  \end{picture}
  \caption {Geometric realizations of trees $T_0^k,T_1^k,\dots,T_{k-1}^k$ for $k=4$} \label{Fig:disks}
\end{figure}

The above theorem could be used in the promised proof of Theorem \ref{problem 1}. But, every instance of the infinite diagram $D_{\infty}\of{X_n,g_n,f_n}$ must be quite complicated since it defines a tree-like continuum without the fixed point property, a very complicated object by itself. Fortunately, our proof of \ref{problem 1} requires only finite (but long) version of the diagram. In fact, there is no advantage in using $D_{\infty}\of{X_n,g_n,f_n}$ instead of $D_{l}\of{X_n,g_n,f_n}$ where $l$ is an arbitrary positive integer. In this section, we construct such finite version of the diagram based on a simpler and easier to follow idea; see \ref{t:longdiagramoftrees}.
Our construction is illustrated in Figure \ref{fig:Tn}. For each pair of integers $\mu$ and $\nu$, let $v_\nu^\mu$ denote the point $\of{\mu,\nu}\in\mathbb{R}^2$. Let $k\ge2$ be an integer. For each $n=0,\dots,k-1$, we define a simplicial tree $T_n^k$ in the following way.
\begin{itemize}
    \item $\mathcal{V}\of{T_n}$, the set of vertices of $T_n^k$, consists of points in the following five non-repetitive ordered groups:
    \begin{enumerate}
      \item $v_0^1,v_1^1,\dots,v_n^1,v_{n+1}^0$, (only $v_0^1,v_1^0$ if $n=0$);
      \item $v_0^{-1},v_1^{-1},\dots,v_n^{-1},v_{n+1}^0$, (only $v_0^{-1},v_1^0$ if $n=0$);
      \item $v_{n+1}^0, v_{n+2}^0, \dots, v_{k}^0$, (only $v_k^0$ if $n=k-1$);
      \item $v_{k}^0,v_{k+1}^1,v_{k+2}^1,\dots, v_{k+1+n}^1$, (only $v_k^{0},v_{k+1}^1$ if $n=0$);
      \item $v_{k}^0,v_{k+1}^{-1},v_{k+2}^{-1},\dots, v_{k+1+n}^{-1}$, (only $v_k^{0},v_{k+1}^{-1}$ if $n=0$).
    \end{enumerate}
    Notice that $v_{n+1}^0$ belongs to the first three groups, $v_{k}^0$ belongs to the last three groups, and the five groups have no other intersections.
    \item $\mathcal{E}\of{T_n}$, the set of edges of $T_n^k$, consists of all edges in the form $\{p,q\}\in \mathcal E\of{T_n^k}$ where $p$ and $q$ are consecutive points in any of the above five groups of vertices.
\end{itemize}

\begin{observation}
 For each $n=0,\dots,k-1$,  $T_n^k$ is a simplicial tree with four endpoints $v_{0}^{1}$, $v_{0}^{-1}$, $v_{k+1+n}^{1}$ and $v_{k+1+n}^{-1}$. If $n<k-1$ then $\ofa{T_n^k}$ is homeomorphic to the letter ``H'' with points of order three at $v_{n+1}^{0}$ and $v_{k}^{0}$. Finally, if $n=k-1$ then $\ofa{T_n^k}$ is homeomorphic to the letter ``X'' with $v_{k}^{0}$ being the point of order four.
\end{observation}

For any $n=0,\dots,k-1$, we define a function $s_{n}:\mathcal{V}\of{T_{n}^k}\to\mathcal{V}\of{T_{n}^k}$ by setting $s_n\of{v_\mu^\nu}=v_\mu^{-\nu}$ for all $v_\mu^\nu\in \mathcal{V}\of{T_{n}^k}$.

\begin{observation}
The function  $s_n$ is a simplicial involution of $T_n^k$ onto itself.
\end{observation}

  For all $n=0,\dots,k-2$, we define functions $\sigma_{n}, \tau_n:\mathcal{V}\of{T_{n+1}^k}\to\mathcal{V}\of{T_{n}^k}$ in the following way. For any $v_\mu^\nu\in \mathcal{V}\of{T_{n+1}^k}$ set
  \begin{equation*}
    \sigma_n\of{v_\mu^\nu}=\left\{
                             \begin{array}{ll}
                               v_{n+1}^0, & \hbox{if $\mu=n+1$;} \\
                               v_{\min\of{\mu,k+n+1}}^\nu, & \hbox{otherwise,}
                             \end{array}
                           \right.
  \end{equation*}
and
  \begin{equation*}
    \tau_n\of{v_\mu^\nu}=\left\{
                             \begin{array}{ll}
                               v_{k}^0, & \hbox{if $\mu=k+1$;} \\
                               v_{\max\of{\mu-1,0}}^\nu, & \hbox{otherwise.}
                             \end{array}
                           \right.
  \end{equation*}

We may summarize the above definitions in the following way. The function $\sigma_n$ collapses each of the two upper edges ($\{v_{k+2+n}^{1},v_{k+1+n}^{1}\}$ and $\{v_{k+2+n}^{-1},v_{k+1+n}^{-1}\}$) of $T_{n+1}^k$ to its lower endpoint ($v_{k+1+n}^{1}$ and $v_{k+1+n}^{-1}$), and moves $v_{n+1}^{1}$ and $v_{n+1}^{-1}$ to $v_{n+1}^{0}$. The function $\tau_n$ leaves the two lowest points $v_{0}^{1}$ and $v_{0}^{-1}$ fixed, and moves each other vertex of $T_{n+1}^k$ one unit down and keeping the image in the same horizontal position except the vertices $v_{k+1}^{1}$ and $v_{k+1}^{-1}$ which are mapped to $v_{k}^{0}$.

\begin{observation}
  For all $n=0,1,2,\dots,k-2$, both functions $\sigma_{n},\tau_n:\mathcal V\of{T_{n+1}^k}\rightarrow \mathcal V\of{T_n^k}$ are simplicial surjections.
\end{observation}

Let $p_y$ denote the projection of $\mathbb{R}^2$ onto the $y$-axis.

\begin{observation}
  $v_0^1$, $v_0^{-1}$, $v_{k+1+n}^{1}$ and $v_{k+1+n}^{-1}$ are the only coincidence points of $p_y\circ\sigma_n$ and $p_y\circ\tau_n$.
The same points are the only coincidence points of $\sigma_n$ and $\tau_n$.
\end{observation}

\begin{observation}
  For all $n=0,\dots,k-2$, the maps $\sigma_{n}$ and $s_n\circ\tau_n$ have no coincidence points on $T_{n+1}^k$.
\end{observation}

\begin{observation}
 $s_n\circ\sigma_{n}=\sigma_{n}\circ s_{n+1}$ and $s_n\circ\tau_{n}=\tau_{n}\circ s_{n+1}$ for all $n=0,\dots,k-2$.
\end{observation}

\begin{observation}
 $\sigma_n\circ\tau_{n+1}=\tau_n\circ\sigma_{n+1}$ for all $n=0,\dots,k-3$.
\end{observation}

Set $\omega_n=s_n\circ\tau_n$. We summarize the above construction with the following theorem.

\begin{proposition}\label{t:longdiagramoftrees}
  For each  integer $k\ge2$ there exists a commutative simplicial diagram $D_{k-1}\of{\mathcal V\of{T_n^k},\sigma_n,\omega_n}$ with no coincidence points such that
\begin{enumerate}
\item for each $n=0,\dots,k-1$, $T_n^k$ is a simplicial tree such that
\begin{enumerate}
\item $\ofa{T_{k-1}^k}$ is homeomorphic to the letter ``X",
\item $\ofa{T_n^k}$ is homeomorphic to the letter ``H" if $n\le k-2$, and
\end{enumerate}
\item the functions $\sigma_n$ and $\omega_n$ are simplicial surjections for all $n=0,\dots,k-2$.
\end{enumerate}
\end{proposition}

\section{Open covers generated by simplicial diagrams}\label{s6}

In this section we present connections between commutative simplicial diagrams $D_{l}\of{\mathcal V\of{G_n},g_n,f_n}$ with no coincidence points and sequences $(\mathcal{U}_n)$ of taut tree-chains following a fixed-point-free pattern.

Suppose $G$ is an arbitrary simplicial graph in $\mathbb R^3$,  and  $\varepsilon$ is a real number such that $0<\varepsilon<1$. For any vertex $v$ of $G$ and any edge $e=\{v,u\}$, let $\operatorname{st}_\varepsilon\of{v,e}$ be the set of points in the form $(1-t)v+tu$ where $0\le t<\varepsilon$.
Let $\operatorname{st}_\varepsilon\of{v,G}$ be the union of $\ofc{v}$ and of all $\operatorname{st} _\varepsilon\of{v,e}$ where $v\in e\in\mathcal{E}\of{G}$. Clearly, $\operatorname{st}_\varepsilon\of{v,G}$ is an open neighborhood of $v$ in $G$.

\begin{observation}\label{o:stetaepsi}
  Let $G$ be a simplicial graph, let $v\in\mathcal{V}\of{G}$ and let $\eta$ and $\varepsilon$ be real numbers such that $1>\eta>\varepsilon>0$. Then
$\operatorname{cl}\of{\operatorname{st}_\varepsilon\of{v,G}}\subset\operatorname{st}_\eta\of{v,G}$.
\end{observation}

\begin{observation}\label{o:st}
  Let $G$ be a simplicial graph, and let $\varepsilon$ be real number such that $1>\varepsilon>0.5$. Then the following statements are true.
  \begin{enumerate}
    \item\label{o:st:winstu} For all $u,w\in\mathcal{V}\of{G}$,
    $$
    w\in\operatorname{st}_\varepsilon\of{u,G} \Longleftrightarrow w=u.
    $$
    \item\label{o:st:intersect2} The following three conditions are equivalent for all $u,v\in\mathcal{V}\of{G}$.
    \begin{itemize}
      \item $\operatorname{st}_\varepsilon\of{u,G}\cap\operatorname{st}_\varepsilon\of{v,G}\ne\emptyset$,
      \item $\operatorname{cl}\of{\operatorname{st}_\varepsilon\of{u,G}}\cap\operatorname{cl}\of{\operatorname{st}_\varepsilon\of{v,G}}\ne\emptyset$, and
      \item $u$ and $v$ are $1$-close in $G$.
    \end{itemize}
    \item\label{o:st:intersect3} for any three distinct vertices $u,v,w\in\mathcal{V}\of{G}$,
     $$
     \operatorname{st}_\varepsilon\of{u,G}\cap\operatorname{st}_\varepsilon\of{v,G}\cap\operatorname{st}_\varepsilon\of{w,G}=\emptyset.
     $$
    \item $\ofi{u,v}\subset\operatorname{st}_\varepsilon\of{u,G}\cup\operatorname{st}_\varepsilon\of{v,G}$ for all $\{u,v\}\in\mathcal{E}\of{G}$.
    \item \label{o:st:cover} Let $\mathcal{W}=\ofc{\operatorname{st}_\varepsilon\of{v,G} \mid v\in\mathcal{V}\of{G}}$. Then $\mathcal W$ is an open cover of $\ofa G$. The nerve of $\mathcal{W}$ is isomorphic to the simplicial graph $G$ with isomorphism established by $\operatorname{st}_\varepsilon\of{v,G}\hookrightarrow v$.
  \end{enumerate}
\end{observation}

\begin{observation}\label{o:GtoHst}
Let $g:\mathcal V\of G\to \mathcal V\of H$ be a simplicial map between simplicial graphs, and let $\varepsilon$ be real number such that $1>\varepsilon>0.5$. Then the following statements are true.
  \begin{enumerate}
    \item\label{o:GtoHst:union} For all $u\in\mathcal{V}\of{H}$, $\ofa{g}^{-1}\of{\operatorname{st}_\varepsilon\of{u,H}}$ is the union of all $\operatorname{st}_\varepsilon\of{w,G}$ such that $g\of{w}=u$.
    \item\label{o:GtoHst:intersect} For all $u,v\in\mathcal{V}\of{H}$ such that $\ofa{g}^{-1}\of{\operatorname{st}_\varepsilon\of{u,H}}\cap \ofa{g}^{-1}\of{\operatorname{st}_\varepsilon\of{v,H}}\ne\emptyset$, there exist vertices $x,y\in\mathcal{V}\of{G}$ such that $g\of{x}=u$, $g\of{y}=v$, and $x$ and $y$ are $1$-close in $G$. Consequently, $u$ and $v$ are $1$-close in $H$.
   \end{enumerate}
\end{observation}

\begin{construction}\label{c:Un}
    Suppose  $l$ is a positive integer and $T_0,T_1,\dots,T_l$ are simplicial trees. For each $n=0,1,\dots,l-1$, let $f_n:\mathcal V\of{T_{n+1}}\to \mathcal V\of{T_n}$ be a simplicial map and let  $g_n:\mathcal V\of{T_{n+1}}\to \mathcal V\of{T_n}$ be a simplicial surjection such that the diagram $D_l\of{\mathcal V\of{T_n},g_n,f_n}$ is commutative, and $f_0$ and $g_0$ have no proximity vertices.

Let $\varepsilon_0,\varepsilon_1,\dots,\varepsilon_l$ be a strictly decreasing sequence of real numbers such that $1>\varepsilon_0$ and $\varepsilon_l>0.5$.
  For all $n=0,1,\dots,l$ and $v\in\mathcal{V}\of{T_n}$, set
  $$
  U_n^v=\ofa{g_{nl}}^{-1}\of{\operatorname{st}_{\varepsilon_n}\of{v,T_n}}.
  $$
   For all $n=0,1,\dots,l$, set
   $$
   \mathcal{U}_n=\ofc{U_n^v\mid v\in\mathcal{V}\of{T_n}}.
   $$
  Finally, for each $n=0,1,\dots,l-1$, let
  $\varphi_n:\mathcal{U}_{n+1}\to\mathcal{U}_n$ be defined by
  $$
  \varphi_n\of{U_{n+1}^v}=U_n^{f_n\of{v}}.
  $$
\end{construction}

\begin{theorem}\label{t:covers}
Assume the notation from Construction \ref{c:Un}. Then the following statements are true.
  \begin{enumerate}
     \item\label{t:covers:clUn+1} $\operatorname{cl}\of{U_{j}^w}\subseteq U_{n}^{g_{nj}\of{w}}$ for all $j=1,\dots,l$, all $n=0,\dots,j-1$, and $w\in\mathcal{V}\of{T_j}$.
    \item\label{t:covers:uinUnv} $w\in U_n^v \Leftrightarrow g_{nl}\of{w}=v$ for all $n=0,\dots,l$, $v\in\mathcal{V}\of{T_n}$ and $w\in\mathcal{V}\of{T_l}$.
    \item\label{t:covers:3eq} Consider the following four conditions:
        \begin{enumerate}
          \item\label{t:covers:3eq1} $U_j^u\cap U_n^v\ne\emptyset$.
          \item\label{t:covers:3eq2} $\operatorname{cl}\of{U_j^u}\cap \operatorname{cl}\of{U_n^v}\ne\emptyset$.
          \item\label{t:covers:3eq2+} There are vertices $w\in g_{jl}{}^{-1}\of{u}$ and $z\in g_{nl}{}^{-1}\of{v}$  $1$-close in $T_l$.
          \item\label{t:covers:3eq3} $g_{nj}\of{u}$ and $v$ are $1$-close in $T_n$.
        \end{enumerate}
Then conditions (\ref{t:covers:3eq1}), (\ref{t:covers:3eq2}), and (\ref{t:covers:3eq2+}) are equivalent and they imply (\ref{t:covers:3eq3}) for all $j=0,\dots,l$, $n=0,\dots,j$, $u\in\mathcal{V}\of{T_j}$ and $v\in\mathcal{V}\of{T_n}$. All four conditions are equivalent if $j=n=0,\dots,l$,  $u\in\mathcal{V}\of{T_n}$ and $v\in\mathcal{V}\of{T_n}$.
    \item\label{t:covers:3vert} $U_n^u\cap U_n^v\cap U_n^w=\emptyset$ for for any three distinct vertices $u,v,w\in\mathcal{V}\of{T_n}$.
    \item\label{t:covers:nerve} For all $n=0,\dots,l$, $\mathcal{U}_n$ is an open cover of $\ofa{T_l}$. The nerve of $\mathcal{U}_n$ is isomorphic to $T_n$  with the isomorphism established by $U_n^u\hookrightarrow u$.
    \item\label{t:covers:phi} Suppose that $n=0,\dots,l-1$, $u\in\mathcal{V}\of{T_n}$ and $v\in\mathcal{V}\of{T_{n+1}}$ are such that $U_{n}^u\cap U_{n+1}^v\ne\emptyset$. Then $\varphi_n\of{ U_{n+1}^v}\cap U_{n}^u=\emptyset$.
    \item\label{t:covers:phi2} Suppose that $j=2,\dots,l-1$, $n=1,\dots,j$,  $u\in\mathcal{V}\of{T_{j+1}}$ and $v\in\mathcal{V}\of{T_{n+1}}$. Then the following two implications are true.
        \begin{enumerate}
          \item\label{t:covers:phi2a} $U_{j+1}^u\cap U_{n+1}^v\ne\emptyset$ $\Rightarrow$ $\varphi_j\of{U_{j+1}^u}\cap \varphi_j\of{U_{n+1}^v}\ne\emptyset$.
          \item\label{t:covers:phi2b} $U_{j+1}^u\subseteq U_{n+1}^v$ $\Rightarrow$ $\operatorname{cl}\of{\varphi_j\of{U_{j+1}^u}}\subseteq \varphi_j\of{U_{n+1}^v}$.
        \end{enumerate}
        \item\label{t:covers:phi3} Suppose that $n=0,\dots,l-1$ and $u,v\in\mathcal{V}\of{T_{n+1}}$ are such that $U_{n+1}^u\cap U_{n+1}^v\ne\emptyset$. Then $\varphi_{n}\of{U_{n+1}^u}\cap\varphi_{n}\of{U_{n+1}^v}\ne\emptyset$.

  \end{enumerate}
\end{theorem}
\begin{proof}

\noindent
Proof of (\ref{t:covers:clUn+1}). By  \ref{o:GtoHst}(\ref{o:GtoHst:union}), $\operatorname{cl}\of{U_{j}^w}$ is the union of all $\operatorname{cl}\of{\operatorname{st}_{\varepsilon_{j}}\of{u,T_l}}$ such that $g_{jl}\of{u}=w$. For each such $u$, $g_{nl}\of{u}=g_{nj}\of{w}$ and $\operatorname{cl}\of{\operatorname{st}_{\varepsilon_{j}}\of{u,T_l}}\subseteq \operatorname{st}_{\varepsilon_{n}}\of{u,T_l}$ by \ref{o:stetaepsi}. Using \ref{o:GtoHst}(\ref{o:GtoHst:union}) again we get the result that $\operatorname{st}_{\varepsilon_{n}}\of{u,T_l}\subseteq U_n^{{g_{nj}}\of{w}}$. Consequently,  $\operatorname{cl}\of{U_{j}^w}\subseteq U_{n}^{g_{nj}\of{w}}$.

\noindent
Proof of (\ref{t:covers:uinUnv}). If $g_{nl}\of{w}=v$, then \ref{o:GtoHst}(\ref{o:GtoHst:union}) implies that $w\in\operatorname{st}_{\varepsilon_{n}}\of{w,T_l}\subseteq U_n^v$.
On the other hand, suppose $w\in U_n^v$. Then \ref{o:GtoHst}(\ref{o:GtoHst:union}) implies that there is $u$ such that $g_{nl}\of{u}=v$ and $w\in\operatorname{st}_{\varepsilon_{n}}\of{u,T_l}$.
Using \ref{o:st}(\ref{o:st:winstu}) we infer that $w=u$ and, consequently, $g_{nl}\of{w}=v$.

\noindent
Proof of (\ref{t:covers:3eq}). The implication (\ref{t:covers:3eq1}) $\Rightarrow$ (\ref{t:covers:3eq2}) is trivial.

\noindent
Suppose (\ref{t:covers:3eq2}). It follows from \ref{o:GtoHst}(\ref{o:GtoHst:union}) that there are $w\in g_{jl}{}^{-1}\of{u}$ and $z\in g_{nl}{}^{-1}\of{v}$ such that
$\operatorname{st}_{\varepsilon_{j}}\of{w,T_l}\subseteq  U_j^u$, $\operatorname{st}_{\varepsilon_{n}}\of{z,T_l}\subseteq U_n^v$ and
$\operatorname{cl}\of{\operatorname{st}_{\varepsilon_{j}}\of{w,T_l}}\cap\operatorname{cl}\of{\operatorname{st}_{\varepsilon_{n}}\of{z,T_l}}\ne\emptyset$. So, $\operatorname{cl}\of{\operatorname{st}_{\varepsilon_{n}}\of{w,T_l}}\cap\operatorname{cl}\of{\operatorname{st}_{\varepsilon_{n}}\of{z,T_l}}\ne\emptyset$ because $\varepsilon_n\ge\varepsilon_j$.
Using \ref{o:st}(\ref{o:st:intersect2}) with $\varepsilon=\varepsilon_n$ we infer that $w$ and $z$ are $1$-close in $T_l$. So, the implication (\ref{t:covers:3eq2}) $\Rightarrow$ (\ref{t:covers:3eq2+}) is true.

\noindent
Suppose (\ref{t:covers:3eq2+}). Since $w$ and $z$ are $1$-close in $T_l$, \ref{o:st}(\ref{o:st:intersect2}) used with $\varepsilon=\varepsilon_j$ implies that $\operatorname{st}_{\varepsilon_{j}}\of{w,T_l}\cap\operatorname{st}_{\varepsilon_{j}}\of{z,T_l}\ne\emptyset$.
Also $\operatorname{st}_{\varepsilon_{j}}\of{w,T_l}\cap\operatorname{st}_{\varepsilon_{n}}\of{z,T_l}\ne\emptyset$  because $\varepsilon_n\ge\varepsilon_j$. Consequently, $U_j^u\cap U_n^v\ne\emptyset$ and the implication (\ref{t:covers:3eq2+}) $\Rightarrow$ (\ref{t:covers:3eq1}) is also true, and (\ref{t:covers:3eq1}), (\ref{t:covers:3eq2}), and (\ref{t:covers:3eq2+}) are equivalent.

\noindent
Still assuming (\ref{t:covers:3eq2+}), we infer $g_{nl}\of{w}=g_{nl}\circ g_{jl}\of{w}=g_{nj}\of{u}$ and $v=g_{nl}\of{z}$ are $1$-close in $T_n$. So, the implication (\ref{t:covers:3eq2+}) $\Rightarrow$ (\ref{t:covers:3eq3}) is true.

\noindent
Proof of (\ref{t:covers:3eq3}) $\Rightarrow$ (\ref{t:covers:3eq1}) for $j=n$.
If $u=v$ then (\ref{t:covers:3eq1}) follows from (\ref{t:covers:uinUnv}). So, we may assume that $u$ and $v$ are adjacent in $T_n$. Since $g_{nl}$ is a simplicial surjection, $T_n$ is a tree and $T_l$ is connected, Observation \ref{o:simpsurtree} implies that there is an edge $\{a,b\}\in \mathcal V\of{T_l}$ such that $g_{nl}\of{a}=u$ and $g_{nl}\of{b}=v$.
Observation \ref{o:GtoHst}(\ref{o:GtoHst:union}) implies that $\operatorname{st}_{\varepsilon_n}\of{a,T_l}\subseteq U_n^u$ and $\operatorname{st}_{\varepsilon_n}\of{b,T_l}\subseteq U_n^v$. Using \ref{o:st}(\ref{o:st:intersect2}) we infer that $\operatorname{st}_{\varepsilon_n}\of{a,T_l}\cap \operatorname{st}_{\varepsilon_n}\of{b,T_l}\ne\emptyset$. Therefore, $U_n^u\cap U_n^v\ne\emptyset$.

\noindent
Proof of (\ref{t:covers:3vert}).  (\ref{t:covers:3vert}) follows from (\ref{t:covers:3eq}), since $T_n$ is a simplicial tree and, consequently, there are no three vertices $u, v, w\in T_n$ such that any two of them are adjacent.

\noindent
Proof of (\ref{t:covers:nerve}). Observe that $U_l^w=\operatorname{st}_{\varepsilon_{l}}\of{w,T_l}$ for each vertex $w$ of $T_l$. So, $\mathcal{U}_l$ is an open cover of $\ofa{T_l}$ by \ref{o:st}(\ref{o:st:cover}).
Using (\ref{t:covers:clUn+1}) we observe that $U_l^w\subseteq U_n^{g_{nl}\of{w}}$ for each vertex $w$ of $T_l$. So, $\mathcal{U}_l$ is a refinement of $\mathcal{U}_n$ for all $n=0,\dots,l$.
Since all sets $U_n^v$ are open by \ref{o:GtoHst}(\ref{o:GtoHst:union}), $\mathcal{U}_n$ is an open cover of $\ofa{T_l}$. The dimension of the nerve of $\mathcal{U}_n$ is $\le1$ by (\ref{t:covers:3vert}). So, the second part of (\ref{t:covers:nerve}) follows from the equivalence of (\ref{t:covers:3eq1}) and (\ref{t:covers:3eq3}).

\noindent
Proof of (\ref{t:covers:phi}). Using (\ref{t:covers:3eq}), we infer that
\begin{equation*}\label{e:uandgn(v)close}
  \text{$u$ and $g_n\of{v}$ are $1$-close in $T_n$.} \tag{$*$}
\end{equation*}
Suppose $\varphi_n\of{ U_{n+1}^v}=U_n^{f_n\of{v}}$ intersects $U_{n}^u$. Using (\ref{t:covers:3eq}) again, we infer that
\begin{equation*}\label{e:uandfn(v)close}
  \text{$u$ and $f_n\of{v}$ are $1$-close in $T_n$.} \tag{$**$}
\end{equation*}
Now, \eqref{e:uandgn(v)close} and \eqref{e:uandfn(v)close} together imply that $f_n\of{v}$ and $g_n\of{v}$ are $2$-close in $T_n$, making $v$ a proximity vertex for $f_n$ and $g_n$. But, there is no proximity vertex for $f_n$ and $g_n$ by \ref{o:nopoximityforfnandgn}. This contradiction completes the proof of (\ref{t:covers:phi}).

\noindent
Proof of (\ref{t:covers:phi2a}). Suppose $U_{j+1}^u\cap U_{n+1}^v\ne\emptyset$. By (\ref{t:covers:3eq2+}), there are vertices $w\in g_{j+1,l}{}^{-1}\of{u}$ and $z\in g_{nl}{}^{-1}\of{v}$  $1$-close in $T_l$.
So, $u=g_{j+1,l}\of{w}$ and $x=g_{j+1,l}\of{z}$ are $1$-close in $T_{j+1}$. Clearly, $g_{n+1,j+1}\of{x}=g_{n+1,l}\of{z}=v$. Observe that  $f_n\of{v}=f_n\circ g_{n+1,j+1}\of{x}=g_{nj}\circ f_j\of{x}$ because the diagram $D_j\of{\mathcal V\of{T_j},g_j,f_j}$ is commutative. The vertices $f_j\of{x}$ and $f_j\of{u}$ are $1$-close in $T_j$ since $x$ and $u$ are $1$-close in $T_{j+1}$ and $f_j$ is simplicial. It follows from \ref{o:simpsurtree} that  there are $1$-close vertices $a,b\in\mathcal{V}\of{T_l}$ such that $g_{jl}\of{a}=f_j\of{u}$ and $g_{jl}\of{b}=f_j\of{x}$.
Clearly, $a\in U_j^{f_j\of{u}}=\varphi_j\of{U_{j+1}^u}$.
Since $g_{nl}\of{b}=g_{nj}\circ f_j\of{x}=f_n\circ g_{n+1,j+1}\of{x}=f_n\of{v}$, we get the result that $b\in U_n^{f_n\of{v}}=\varphi_n\of{U_{n+1}^v}$.
Now, (\ref{t:covers:3eq}) implies that $\varphi_j\of{U_{j+1}^u}\cap \varphi_j\of{U_{n+1}^v}\ne\emptyset$.

\noindent
Proof of (\ref{t:covers:phi2b}). Take a vertex $t\in\mathcal{V}\of{T_l}$ such that $g_{j+1,l}\of{t}=u$.
Then $t\in U_{j+1}^u$ by (\ref{t:covers:uinUnv}).
So, $t\in U_{n+1}^v$ and $g_{n+1,l}\of{t}=v$  by (\ref{t:covers:uinUnv}) again.
Observe that $g_{n+1,j+1}\of{u}=v$ because $g_{n+1,j+1}\of{u}=g_{n+1,j+1}\circ g_{j+1,l}\of{t}=g_{n+1,l}\of{t}=v$.
Consequently, $\operatorname{cl}\of{\varphi_{j}\of{U_{j+1}^u}}=\operatorname{cl}\of{U_{j}^{f_{j}\of{u}}}\subseteq U_n^{g_{nj}\circ f_{j}\of{u}}$ by (\ref{t:covers:clUn+1}).
Since $$g_{nj}\circ f_{j}\of{u}=f_n\circ g_{n+1,j+1}\of{u}=f_n\of{v},$$ we infer that $\operatorname{cl}\of{\varphi_{n+1}\of{U_{n+1}^u}}\subseteq U_n^{f_n\of{v}}=\varphi_{n}\of{U_{n+1}^v}$.

\noindent
Proof of (\ref{t:covers:phi3}). (\ref{t:covers:3eq}) implies that $u$ and $v$ are $1$-close in $T_{n+1}$. So, $f_n\of{u}$ and $f_n\of{v}$ are $1$-close in $T_{n}$. Now, we use implication (\ref{t:covers:3eq3}) $\Rightarrow$ (\ref{t:covers:3eq1}) for $j=n$ to infer that $U_n^{f_n\of{u}}\cap U_n^{f_n\of{v}}\ne\emptyset$. Consequently, $\varphi_{n}\of{U_{n+1}^u}\cap\varphi_{n}\of{U_{n+1}^v}\ne\emptyset$.
\end{proof}

\begin{theorem}\label{t:coversofX}
  Let $l\ge 2$ be an integer and let $X$ be a topological tree homeomorphic to the letter ``X.'' Then there exists a sequence $\mathcal{U}_0,\mathcal{U}_1,\dots,\mathcal{U}_l$ of finite tree-covers of $X$, and there exists a sequence of functions $\varphi_0:\mathcal{U}_1\to\mathcal{U}_0, \varphi_1:\mathcal{U}_2\to\mathcal{U}_1,\dots,\varphi_{l-1}:\mathcal{U}_l\to\mathcal{U}_{l-1}$ satisfying the following conditions.
   \begin{enumerate}
  \item\label{t:coversofX1} For each $n=0,\dots,l-1$, $\mathcal{U}_{n+1}$ strongly refines $\mathcal{U}_{n}$, i.e. for each $U\in\mathcal{U}_{n+1}$, there is $V\in\mathcal{U}_n$ such that $\operatorname{cl}\of{U}\subseteq V$.
  \item\label{t:coversofX2}  For each $n=0,\dots,l-1$, for each $U\in\mathcal{U}_n$ and for each $V\in\mathcal{U}_{n+1}$,
  $$U\cap V\ne\emptyset \quad \Rightarrow \quad  \varphi_n(V)\cap U=\emptyset.$$
   \item\label{t:coversofX3} Suppose that $j=2,\dots,l-1$, $n=1,\dots,j$,  $U\in\mathcal{U}_{j+1}$ and $V\in\mathcal{U}_{n+1}$. Then
  \begin{enumerate}
    \item\label{t:coversofX3a} $U\cap V\ne\emptyset \quad \Rightarrow \quad  \varphi_{n+1}(U)\cap \varphi_{n}(V)\ne\emptyset$, and
    \item\label{t:coversofX3b} $\operatorname{cl}(U)\subseteq V\ \quad \Rightarrow \quad  \varphi_{n+1}(U)\subseteq \varphi_{n}(V)$.
  \end{enumerate}
  \item\label{t:coversofX4} For each $n=0,\dots,l-1$, and for all $U,V\in\mathcal{U}_{n+1}$,
   $$U\cap V\ne\emptyset \quad \Rightarrow \quad  \varphi_{n}(U)\cap \varphi_n(V)\ne\emptyset.$$
  \item\label{t:coversofX5} The collection $\mathcal{U}=\bigcup_{n=0}^l\mathcal{U}_n$ is taut.
\end{enumerate}
\end{theorem}
\begin{proof}
 Use Proposition \ref{t:longdiagramoftrees} with $k=l+1$ to get a simplicial, surjective and commutative diagram $D_l\of{\mathcal V\of{Y_n},p_n,q_n}$ such that $\ofa{p_0}$ and $\ofa{q_0}$  have no coincidence points, $Y_l=T_l^{l+1}$ and its geometric realization $\ofa{T_l^{l+1}}$ is homeomorphic to the letter ``X", and $Y_n=T_n^{l+1}$ and its geometric realization $\ofa{T_n^{l+1}}$ is homeomorphic to the letter ``H" for all $n=0,\dots,l-1$. Now, set $D_l\of{\mathcal V\of{T_n},g_n,f_n}=D_l\of{\mathcal V\of{Y_n^{\of{3}}},p_n^{\of{3}},q_n^{\of{3}}}$ and  use Corollary \ref{t:nopoximitydiagram} to infer that $D_l\of{\mathcal V\of{T_n},g_n,f_n}$
is a simplicial, surjective and commutative diagram  such that
\begin{itemize}
  \item $g_0$ and $f_0$ have no proximity vertices,
  \item $\ofa{T_n}$ is homeomorphic to the letter ``H" for all $n=0,\dots,l-1$, and
  \item $\ofa{T_l}$ is homeomorphic to the letter ``X."
\end{itemize}
We can  now use \ref{c:Un} to construct $\mathcal{U}_0,\mathcal{U}_1,\dots,\mathcal{U}_l$ and $\varphi_0, \varphi_1,\dots,\varphi_{l-1}$ satisfying the properties listed in Theorem \ref{t:covers}.

Observe that \ref{t:covers} (\ref{t:covers:nerve}) implies that all $\mathcal{U}_0,\mathcal{U}_1,\dots,\mathcal{U}_l$ are finite tree-covers of $X$.
To complete the proof of Theorem \ref{t:coversofX} we need to observe that (\ref{t:coversofX1})-(\ref{t:coversofX5}) are  satisfied.
(\ref{t:coversofX1}) follows from \ref{t:covers} (\ref{t:covers:clUn+1}).
(\ref{t:coversofX2}) follows from \ref{t:covers} (\ref{t:covers:phi}). (\ref{t:coversofX3})  follows from \ref{t:covers} (\ref{t:covers:phi2}). (\ref{t:coversofX4})  follows from \ref{t:covers} (\ref{t:covers:phi3}).
Finally,  (\ref{t:coversofX5})  follows from the equivalence of \ref{t:covers} (\ref{t:covers:3eq1}) and \ref{t:covers} (\ref{t:covers:3eq2}).
\end{proof}

\begin{proof}[Proof of Theorem \ref{problem 1}]
 To complete the proof place a continuum $X$ homeomorphic to the letter ``X" in the plane, use Theorem \ref{t:coversofX} to get a collection $\mathcal{U}$ promised by the theorem. Finally, use Proposition \ref{p:enlarge} with $M=\mathbb{R}^2$ to replace each $U\in\mathcal{U}$ by $\widetilde U$.
 \end{proof}

\begin{remark}
  Observe that Proposition \ref{t:longdiagramoftrees} can be replaced in the proof of \ref{problem 1} by any  commutative, surjective and long simplicial diagram of trees with no coincidence points. So, we could use the infinite diagrams ($D_{\infty}\of{X_n,g_n,f_n}$ where $X_n$'s are simplicial trees) constructed by  Oversteegen and Rogers in \cite[Section 3]{39}, and by Hernández-Gutiérrez and  Hoehn in \cite{logan}; see Theorem \ref{t:infdiagramoftrees}. However, using an infinite diagram gives no apparent advantage because of the following reason.
  In the above construction, an embedding $h_l:X_l\to\mathbb{R}^2$ is selected and then covers $\mathcal{U}_0^l,\mathcal{U}_1^l,\dots,\mathcal{U}_l^l$ are constructed heavily depending on $h_l$. In order to extend this construction, one would have to take an integer $m>l$, select an embedding $h_m:X_l\to\mathbb{R}^2$ and then obtain covers 
  $\mathcal{U}_0^m,\mathcal{U}_1^m,\dots,\mathcal{U}_m^m$. However, there is no apparent nexus between the two sequences of covers because there is no connection between the embeddings $h_l$ and $h_m$. In particular, since $\varprojlim\of{X_n,g_n}$ is not embeddable in the plane, there may be no embedding $h_m:X_l\to\mathbb{R}^2$ such that $h_m$ is ``sufficiently close'' to $h_l\circ g_{lm}$.

\end{remark}


\noindent I. Bani\v c\\
              (1) Faculty of Natural Sciences and Mathematics, University of Maribor, Koro\v{s}ka 160, SI-2000 Maribor,
   Slovenia; \\(2) Institute of Mathematics, Physics and Mechanics, Jadranska 19, SI-1000 Ljubljana, 
   Slovenia; \\(3) Andrej Maru\v si\v c Institute, University of Primorska, Muzejski trg 2, SI-6000 Koper,
   Slovenia\\
             {iztok.banic@um.si}           
     
				\-

  \noindent J.  Kennedy\\
             Lamar University, 200 Lucas Building, P.O. Box 10047, Beaumont, TX 77710 USA\\
{{kennedy9905@gmail.com}       }    

                 	\-

  \noindent P.  Minc\\
             Department of Mathematics, 218 Parker Hall, Auburn University, Auburn, AL 36849-5310, USA\\
{{mincpio@auburn.edu}       }    

\end{document}